\documentclass[opre,nonblindrev]{myinforms} %

\OneAndAHalfSpacedXI %

\usepackage{algorithm,algpseudocode}

\usepackage{natbib}
\bibpunct[, ]{(}{)}{,}{a}{}{,}%
\def\bibfont{\small}%

\TheoremsNumberedThrough     %
\ECRepeatTheorems

\EquationsNumberedThrough    %

\usepackage[utf8]{inputenc} %
\usepackage[T1]{fontenc}    %
\usepackage{hyperref}       %
\usepackage{url}            %
\usepackage{booktabs}       %
\usepackage{amsfonts}       %
\usepackage{nicefrac}       %
\usepackage{microtype}      %
\usepackage{xcolor}         %
\usepackage{amsmath}
\usepackage{amssymb}
\usepackage{multirow}
\usepackage{makecell}
\usepackage{algorithm}
\usepackage{algpseudocode}
\usepackage[capitalise]{cleveref}
\usepackage{shortcuts}
\usepackage{booktabs}
\usepackage{array, caption, threeparttable}
\usepackage{caption}
\usepackage{tikz}
\usetikzlibrary{positioning,arrows,decorations,decorations.markings,shadows,positioning,arrows.meta,matrix,fit}

\Crefname{assumption}{Assumption}{Assumptions}
\Crefname{remark}{Remark}{Remarks}

\usepackage[compact]{titlesec}

\begin{document}
		
	\RUNAUTHOR{Li et al.}
	
	\RUNTITLE{A Fenchel-Young Loss Approach to Data-Driven Inverse Optimization}
	
	\TITLE{A Fenchel-Young Loss Approach to Data-Driven Inverse Optimization}
	
	\newcommand*\samethanks[1][\value{footnote}]{\footnotemark[#1]}
	\ARTICLEAUTHORS{%
		\AUTHOR{Zhehao Li\thanks{These authors contributed equally to this work.}}
 		\AFF{Tsinghua University, 100084 Beijing, China, \EMAIL{lizh21@mails.tsinghua.edu.cn}}
        \AUTHOR{Yanchen Wu\samethanks}
 		\AFF{Tsinghua University, 100084 Beijing, China, \EMAIL{wu-yc23@mails.tsinghua.edu.cn}}
        \AUTHOR{Xiaojie Mao\footnote{Corresponding author.}}
 		\AFF{Tsinghua University, 100084 Beijing, China, \EMAIL{maoxj@sem.tsinghua.edu.cn}}
} %

\ABSTRACT{\noindent
    Data-driven inverse optimization seeks to estimate unknown parameters in an optimization model from observations of optimization solutions. 
    Many existing methods are ineffective in handling noisy and suboptimal solution observations and also suffer from computational challenges.
    In this paper, we build a connection between inverse optimization and the Fenchel-Young (FY) loss originally designed for structured prediction, proposing a FY loss approach to data-driven inverse optimization. 
    This new approach is amenable to efficient gradient-based optimization, hence much more efficient than existing methods. We provide theoretical guarantees for the proposed method and use extensive simulation and real-data experiments to demonstrate its significant advantage in parameter estimation accuracy, decision error and computational speed. 
}

\maketitle

\section{Introduction}

    Inverse optimization is a fundamental problem in operations research and machine learning, aiming to estimate unknown parameters in an optimization model from observations of solutions \citep{ahuja2001inverse}. Inverse optimization is important in many applications, such as consumer utility modeling \citep{saez2016data,saez2017short}, vehicle routing \citep{chen2021inverse,zattoni2024inverse}, medical decision-making \citep{chan2014generalized,aswani2019data}, and portfolio risk management \citep{li2021inverse,yu2023learning}, etc. 
     For example, in the consumer utility modeling problem,  customers' utility function may involve unknown parameters characterizing their preference for different products. These parameters may be estimated from observations of utility-maximizing customers' purchase decisions.  
    By recovering the unknown parameters in the optimization, analysts can better understand the underlying decision-making process and make better decisions of their interest.


    In data-driven inverse optimization, the observed decisions may be often noisy, suboptimal and even infeasible, which poses significant challenges for estimating the unknown parameters \citep{chan2023inverse}. 
    \citet{aswani2018inverse} discovered that many existing approaches may not effectively recover the true parameter value, such as the Karush-Kuhn-Tucker approximation (KKA) approach  \citep{keshavarz2011imputing} and the variational inequality approximation (VIA) approach  \citep{bertsimas2015data}. \citet{aswani2018inverse} then proposed a parameter estimator such that the average squared distance between the decision induced by the estimated parameter and the observed decision is minimized and showed its effectiveness theoretically and numerically.
    However, this  approach involves non-convex optimization  solved by brute-force enumeration, restricted to only small-scale problems. 
    For larger problems, \citet{aswani2018inverse} relied on a semiparametric approach that is more scalable than enumeration. However, this approach  it is not amenable to first-order gradient-based optimization, hence still computationally challenging for large-scale problems.  
    
    

    In this paper, we propose an alternative surrogate loss function approach to overcome the drawbacks of existing methods. Specifically, we consider a regularized version of the optimization problem of interest and adapt the suboptimality loss in
\citet{chan2014generalized,chan2018trade} to the regularized problem. We find that the resulting loss coincides with the Fenchel-Young (FY) loss originally designed for structured prediction tasks such as classification and ranking \citep{blondel2020learning, blondel2022learning}.
This builds a novel connection between FY loss and data-driven inverse optimization, providing new perspectives on both. Importantly, the FY loss is differentiable with a closed-form gradient and it is convex for linear problems. 
This allows for efficient first-order optimization like stochastic gradient descent (SGD), thereby achieving significant computational advantage over existing methods. 
By extending existing FY loss theory, we also provide some theoretical guarantees for the proposed method. We evaluate our method on synthetic and real-world datasets, demonstrating that our approach outperforms existing methods in terms of  parameter estimation accuracy, generating high quality decisions,  and computational efficiency. In particular, the proposed method can achieve robust performance even in settings where existing methods easily give degenerate solutions. 

\vskip 0.2in

    \subsection{Related Work}

        \textbf{Inverse Optimization} ~~
        Early work on inverse optimization focused on parameter recovery from  noiseless observations of decision solutions \citep{ahuja2001inverse, chan2014generalized, chan2018trade}. \citet{barmann2017emulating} studied inverse optimization with noiseless solutions in the online setting, focusing on minimizing decision suboptimality, while \citet{besbes2023contextual} aimed to estimate parameters by minimizing worst-case regret. To handle noisy and suboptimal decisions in data-driven inverse-optimization,  \citet{aswani2018inverse} proposed a squared distance loss approach and a semiparametric approach and compared their methods with methods in \citet{keshavarz2011imputing,bertsimas2015data}.
         Their methods were also expanded to the online setting \citep{dong2018generalized}. 
        \citet{chan2019inverse} proposed a goodness-of-fit measure for data-driven inverse-optimization.
        We refer to \citet{chan2023inverse} for a comprehensive survey of literature on classical inverse optimization and data-driven inverse optimization and, \citet{mohajerin2018data,birge2022stochastic,lin2024conformal} for some recent advances in enhancing the robustness of inverse optimization methods. Our paper considers the offline setting in \citet{aswani2018inverse} and proposes a Fenchel-Young loss approach to inverse optimization, addressing the computational challenges with existing methods while achieving rigorous theoretical guarantees.

\vskip 0.1in
    
        \noindent \textbf{Structured Prediction with FY Loss} ~~
        The Fenchel-Young (FY) loss was originally introduced for structured prediction (SP) tasks such as multi-class classification and ranking \citep[e.g.,][]{blondel2019learning,blondel2020learning, blondel2022learning}. We refer to \citet{blondel2020learning} for a systematic introduction to FY loss. In this paper, we build the connection between FY loss and data-driven inverse-optimization. Our methodology and theory  build on and extend the FY loss literature, broadening the applicability of FY loss. In particular, the traditional structured prediction tasks for which FY loss is designed  typically involve optimization or ranking over simple sets like simplex, while inverse optimization  may often involve operations problems with more complex structure. 
        

\vskip 0.3in

        \section{Problem Formulation}\label{sec:formulation}

        \subsection{Forward and Inverse Optimization}
    
            In this paper, we study the inverse linear optimization problem with noisy observed decisions. Let $x \in \mathbb{R}^d$ be the decision vector, $u \in \mathbb{R}^m$ the context variable, and $\theta^{\star} \in \mathbb{R}^p$ the unknown parameter vector to be estimated. The forward optimization problem (FOP) is given by:
            \begin{equation}\label{eq:FOP}
                \operatorname{FOP}(\theta^{\star}, u):=  \max_{x} \{h(\theta^{\star}; u)^{\top} x \mid g(u, x) \preceq 0 \}.
            \end{equation}
            Here $h: \mathbb{R}^p \times \mathbb{R}^m \to \mathbb{R}^d$  represents the cost function, and function $g: \mathbb{R}^m \times \mathbb{R}^d \to \mathbb{R}^q$ defines the constraints. Both $h$ and $g$ are assumed known, leaving $\theta^{\star}$ as the only unknown parameter. 
            We define the feasible region of FOP as $\mathcal{X}(u) := \{x \in \mathbb{R}^d \mid g(u, x) \preceq 0 \}$, and use $x^{\star}(\theta^{\star}; u) \in \mathcal{X}^{\star}(\theta^{\star}; u)$  to denote a generic optimal solution  of $\operatorname{FOP}(\theta^{\star}, u)$, where $\mathcal{X}^{\star}(\theta; u) := \argmax_{x}\{h(\theta; u)^{\top} x \mid g(u, x) \preceq 0 \}$ is the optimal solution set.
            
            We assume that $h(\theta; u)$ is continuous in $\theta$ and $u$, while $g(u, x)$ is continuous in both $u$ and $x$ and convex in $x$ for each fixed $u$. These conditions mean that the FOP in \Cref{eq:FOP} has a linear objective with convex constraints\footnote{While we focus on linear objectives in the main text, our method can be also applied to nonlinear objectives; see the numerical experiments in \Cref{sec: Experiment-details-appendix}.}. Such formulation is widely studied in the inverse optimization literature \citep{bertsimas2015data, barmann2017emulating, aswani2018inverse, chen2021inverse, yu2023learning}. Additional examples that align with our formulation in \Cref{eq:FOP} are provided in \Cref{sec:supplementary-formulation}.

\vskip 0.1in
    
            \noindent \textbf{Noisy Observations} ~~ Suppose $(U, Y) \in \mathbb{R}^{m} \times \mathbb{R}^d$ is a pair of context-decision variables, sampled from an unknown joint distribution $P$. Here $Y$ stands for the noisy observation of the true decision $x^{\star}(\theta^{\star}; U)$ under the context $U$ with the unknown parameter $\theta^{\star}$.
            We assume observing a sequence of noisy data $\{(U_i, Y_i)\}_{i=1}^{n}$ drawn identically and independently (IID) from the joint distribution $P$. 
    
            \begin{assumption}[Noisy Decisions]\label{asmp:noisy-decisions}
                For each $i \in [n]$, we assume $Y_i = x^{\star}(\theta^{\star}; U_i) + W_i$. The noise terms $\{W_i\}_{i=1}^{n}$ are IID with zero conditional mean and finite conditional variance given context $u \in \mathbb{R}^m$, i.e., $\mathbb{E}[W_i \mid U = u] = 0$  and $\var[W_i \mid U = u] < \infty$. Meanwhile, $W_i$ is independent of $U_i$ for each $i \in [n]$.
            \end{assumption}
    
            As highlighted in \citet{aswani2018inverse}, accounting for noise in inverse optimization is crucial, as real-world decisions are often influenced by various sources of noise. This noise can arise from (i) measurement errors in the data collection process, (ii) deviations of experts from optimal behavior, or (iii) discrepancies between the parametric FOP model and the true underlying decision-making process.
    
            \begin{remark}[Noise Model]\label{rmk:noisy-objective}
                In this paper, we follow \citet{aswani2018inverse} and  focus on the noisy decisions setting in \Cref{asmp:noisy-decisions}. An alternative noisy model is to incorporate noise in the optimization objective, i.e., we observe $Y_i \in \argmax_{x \in \mathcal{X}(U_i)} \{(h(\theta^{\star}; U_i) + W_i)^{\top} x \}$ for each $i \in [n]$. Our experiments in \Cref{sec: Experiment-details-appendix} 
                also examine this noisy objective setting and demonstrate that our proposed method also outperforms existing methods in this setting.
            \end{remark}
            

\vskip 0.1in
    
            \noindent \textbf{Data-driven Inverse Optimization} ~~ Given a sequence of data $\{(U_i, Y_i)\}_{i=1}^{n}$, the goal of inverse optimization is to estimate the unknown parameter $\theta^{\star}$ from data. The ultimate aim is to leverage the estimated parameter, denoted as $\hat{\theta}$, to make new decisions in subsequent contexts. Therefore, it is desirable to ensure that the optimal decision $x^{\star}(\hat{\theta}; U)$ is, on average, as close as possible to the observed decision $Y$. To achieve this goal, \citet{aswani2018inverse} then proposed to minimize the risk associated with the inverse optimization problem (IOP):
            \begin{equation}\label{eq:IOP-risk}
                \min_{\theta \in \Theta} ~ R_{\rm DIST}(\theta) := \mathbb{E}\big[L_{\rm DIST}(\theta; U, Y) \big],
            \end{equation}
            where $\Theta \in \mathbb{R}^p$ and $L_{\rm DIST}$ is the squared distance loss:
            \begin{equation}\label{eq:distance-loss}
                L_{\rm DIST}(\theta; U, Y) := \min_{x \in \mathcal{X}^{\star}(\theta; U)} \big\|Y - x \big\|_2^2,
            \end{equation} 
            which quantifies the deviation of the optimal solution set $\mathcal{X}^{\star}(\theta; U)$ from the observed decision $Y$ under the context $U$. \Cref{asmp:noisy-decisions} guarantees that the true unknown parameter $\theta^{\star}$ minimizes the risk $R_{\rm DIST}(\theta)$ \citep{aswani2018inverse}.
            As an alternative to the IOP formulation in \Cref{eq:IOP-risk}, we can consider minimizing the excess risk $R_{\rm DIST}(\theta) - R_{\rm DIST}(\theta^{\star})$ over $\theta$. 
            It turns out that the excess risk is equivalent to the expected decision error (defined below) over the distribution of $U$ when the optimal solution of FOP with respect to any given parameter $\theta$ is almost surely unique.
            
            \begin{assumption}[Almost Sure Uniqueness, \citet{aswani2018inverse}, IC condition]\label{asmp:almost-sure-uniqueness}
                For any $\theta \in \Theta$ and $u \in \mathbb{R}^m$,  $\operatorname{FOP}(\theta, u)$ has a unique optimal solution almost surely, that is, $\mathcal{X}^{\star}(\theta; u) = \{x^{\star}(\theta; u) \}$ almost surely.
            \end{assumption}
    
            \Cref{asmp:almost-sure-uniqueness} is a common assumption \citep{aswani2018inverse,elmachtoub2022smart} and holds automatically when the context $U$ follows a continuous distribution and the constraint set $\mathcal{X}(u)$ is a polytope.
    
            \begin{lemma}\label{lm:excess-distance-risk}
                Given \Cref{asmp:noisy-decisions,asmp:almost-sure-uniqueness}, we have
                $$
                    R_{\rm DIST}(\theta) - R_{\rm DIST}(\theta^{\star}) = \mathbb{E}[ \|x^{\star}(\theta; U) - x^{\star}(\theta^{\star}; U)\|_2^2 ], ~~ \forall ~ \theta \in \Theta.
                $$
            \end{lemma}
    
            \Cref{lm:excess-distance-risk} shows that, when choosing the squared distance loss in \Cref{eq:distance-loss}, minimizing the IOP risk $R(\theta)$ over $\Theta$ is equivalent to minimizing the expected \emph{decision error}
            \begin{equation}\label{eq:expected-decision-error}
                D(\theta, \theta^{\star}) := \mathbb{E}[\|x^{\star}(\theta; U) - x^{\star}(\theta^{\star}; U)\|_2^2],
            \end{equation}
            which is the expectation of the squared distance between the optimal solutions  $x^{\star}(\theta; U)$ and $x^{\star}(\theta^{\star}; U)$ over the distribution of $U$. 
            If a parameter $\theta$ can achieve a small expected decision error, then the decision induced by $\theta$ would achieve a small \emph{decision regret} relative to the decision induced by the true parameter $\theta^\star$:             
    \begin{equation}\label{eq:decision-regret}
                \operatorname{Reg}(\theta, \theta^{\star}) := \mathbb{E}[ h(\theta^{\star}; U)^{\top} (x^{\star}(\theta; U) - x^{\star}(\theta^{\star}; U) ) ].
            \end{equation}
            In \Cref{thm:decisioin-regret-bound}, we show that the decision regret above can be upper bounded by the expected decision error. 
            If we further assume a separation condition as formalized in \citet[IC condition]{aswani2018inverse}, 
            then we can also bound the following \emph{parameter error} through the decision error $D(\theta, \theta^{\star})$:
            \begin{equation}\label{eq:parameter-error}
                E(\theta, \theta^{\star}) := \|\theta - \theta^{\star}\|_1.
            \end{equation}
            The parameter error is defined as the distance between a given parameter $\theta$ and the true unknown parameter $\theta^{\star}$ with respect to the $l_q$-norm on the $\mathbb{R}^{p}$ space. Without loss of generality, we choose $q = 1$ for the parameter error $E(\theta, \theta^{\star})$ as shown in \cref{eq:parameter-error}.

\vskip 0.1in

            \noindent \textbf{Empirical Risk Minimization (ERM) with Distance Loss} ~~ \citet{aswani2018inverse} propose to minimize the empirical risk with the squared distance loss in \Cref{eq:distance-loss} over a  parameter space $\Theta$:
            \begin{equation}\label{eq:empirical-risk-minimization}
                \min_{\theta \in \Theta} ~ \widehat{R}_{\rm DIST}(\theta) := \frac{1}{n} \sum_{i=1}^{n} L_{\rm DIST}(\theta; U_i, ~Y_i).
            \end{equation}
            However, directly minimizing the empirical risk in \Cref{eq:empirical-risk-minimization} presents challenges, as the distance loss $L_{\rm DISR}$ is neither convex nor differentiable with respect to the parameter $\theta$, even under \Cref{asmp:almost-sure-uniqueness}. \citet{aswani2018inverse} therefore proposed an enumeration algorithm that computes the value of $\widehat{R}_{\rm DIST}(\theta)$ for each $\theta$ over a discretized grid of $\Theta$, which is computationally expensive when facing the high-dimensional parameter. Meanwhile, the complexity of computing $\widehat{R}_{\rm DIST}(\theta)$   scales significantly with the sample size and decision dimension, further exacerbating computational challenges.
            
            
    
            The next subsection will describe several existing alternatives to the distance loss $L_{\rm DIST}$. To more clearly illustrate properties of these losses, we now assume the cost function $h(\theta; u)$ is linear in the parameter $\theta$. 
            This is the most common setting considered in the existing inverse optimization literature \citep[e.g.,][]{chan2014generalized,barmann2017emulating,chen2021inverse,besbes2023contextual,sun2023maximum}. All numeric examples in  
            \citet{aswani2018inverse} and \Cref{ex:consumer-demand-models,ex:vehicle-routing-problem} also satisfy this linearity condition.

    
            \begin{assumption}[Linear Cost]\label{asmp:linear-cost-function}
                The cost function $h(\theta; u)$ is linear in the unknown parameter $\theta$ for each context $u \in \mathbb{R}^m$.
                , i.e., $h(\theta; u) = A(u) \theta$ for some feature mapping $A(u)$.  
            \end{assumption}

\vskip 0.2in
        
        \subsection{Loss Functions for Inverse Optimization} \label{sub:LossFunction}
    
            The distance loss $L_{\rm DIST}$ is not the only choice for the inverse optimization in \Cref{eq:IOP-risk} via empirical risk minimization. We briefly summarize other loss functions in the existing literature.

\vskip 0.1in
    
            \noindent \textbf{KKT-Approximated Loss} ~~ \citet{keshavarz2011imputing} proposed a loss for a noiseless decision setting that measures the violation of the Karush-Kuhn-Tucker (KKT) conditions:
            $$
            L_{\rm KKA}(\theta; u, y) := \min_{\lambda \succeq 0} \{ L_{\rm ST}(\theta; u, y, \lambda) + L_{\rm CS}(u, y, \lambda) \},
            $$
            where $L_{\rm ST}$ and $L_{\rm CS}$ capture the magnitude of the violation of \emph{stationary condition} and \emph{complementary slackness}, respectively. In noisy setting, minimizing the empirical IOP risk with $L_{\rm KKA}$ is challenging, as observed decisions may be infeasible, leading to large violations of KKT conditions. 

\vskip 0.1in
    
            \noindent \textbf{Suboptimality Loss} ~~ \citet{chan2014generalized,chan2018trade,chan2019inverse} proposed to estimate the parameters by minimizing the average duality gap between the dual optimal objective $v(\theta; \lambda_i, U_i)$ induced by the parameter and the objective function of observed decisions:
            \begin{equation}\label{eq:duality-gap}
                \min_{\theta, \epsilon_i, \lambda_i \geq 0, i= 1, \dots, n } \frac{1}{n} \sum_{i=1}^{n} \epsilon_i ~ \text{s.t.} ~ v(\theta; \lambda_i, U_i) - h(\theta; U_i)^{\top} Y_i \leq \epsilon_i. \\
            \end{equation}
            \citet{chan2023inverse} points out that \Cref{eq:duality-gap} is equivalent to minimizing the empirical risk $\widehat{R}(\theta)$ in \Cref{eq:empirical-risk-minimization} with the loss replaced by the \emph{suboptimality loss} as follows: 
            \begin{equation}\label{eq:Sub-loss}
                L_{\rm SUB}(\theta; u, y) := \max_{x \in \mathcal{X}(u)} h(\theta; u)^{\top} x - h(\theta; u)^{\top} y.
            \end{equation}
            Although the optimization problem in this approach is convex and tractable, this duality-based approach still faces several challenges. (1) If there exists a parameter $\theta \in \Theta$ such that $h(\theta; u) = 0$, then the suboptimality loss $L_{\rm SUB}$ also reaches zero at this $\theta$, causing the inverse optimization problem to degenerate. (2) As the sample size and decision dimension increase, the dual formulation of  \Cref{eq:duality-gap}  becomes computationally challenging. (3) It is unknown whether ERM with the suboptimality loss can control the decision error in \Cref{eq:expected-decision-error}.
            
    

\vskip 0.1in
            
            \noindent \textbf{Variational Inequality Loss} ~~ \citet{bertsimas2015data} study the variational inequality approximation (VIA) loss function that measures the optimality margin for the FOP with a differentiable objective function $f(\theta; u, x)$:
            \begin{equation}\label{eq:VIA-loss}
                L_{\rm VIA}(\theta; u, y) = \max_{x \in \mathcal{X}(u)} \nabla_y f(\theta; u, y)^{\top} (x - y).
            \end{equation}
            If the objective function is linear in decision $x$, e.g., $f(\theta; u, x) = h(\theta; u)^{\top} x$, the VIA loss simplifies to the suboptimality loss, inheriting the same challenges.
            

\vskip 0.1in
    
            \noindent \textbf{Semi-parametric Methods} ~~ To overcome the limitations of squared distance loss with large problems, \citet{aswani2018inverse} propose a semi-parametric (SPA) method that first uses the Nadaraya-Watson (NW) estimator to denoise the observed decisions and then projects them onto the feasible region, denoted as $\tilde{Y}_i$. 
            \citet{aswani2018inverse} then propose to solve \Cref{eq:duality-gap} by replacing $Y_i$ with $\tilde{X}_i$ for each $i \in [n]$.
            However, the NW estimator performs poorly with high-dimensional context space, and they also employ a dual formulation that is computationally challenging to optimize when the sample size and dimension are large.
    
\vskip 0.1in
    
            \noindent \textbf{Challenges for Existing Loss Functions} ~~
            Existing loss functions face several key challenges: difficulty with handling noisy observations, degeneracy, and computational bottleneck. 
            Moreover, most of these methods, except the distance loss and semi-parametric method proposed in \citep{aswani2018inverse}, do not provide direct control over the expected decision error $D(\theta, \theta^{\star})$ in \Cref{eq:expected-decision-error} and parameter error $E(\theta, \theta^{\star})$ in \Cref{eq:parameter-error}. In this paper, we propose a surrogate loss that is computationally efficient, robust to both noise and degeneracy, and retains statistical guarantees for error control.
    
\vskip 0.3in
    
    \section{Inverse Optimization with FY Loss}\label{sec:Fenchel-young-loss}
    
        In this section, we first define the regularized forward optimization problem (R-FOP) by introducing a strongly-convex regularizer, and then connect the suboptimality loss for the R-FOP with the Fenchel-Young (FY) loss. The excess risk induced by FY loss also provides a direct control of the decision error $D(\theta, \theta^{\star})$ for any parameter $\theta$. The convexity and differentiability of FY loss enable us to apply a gradient-based algorithm that is particularly suitable for large sample size and high-dimensional problem. 

\vskip 0.2in
    
        \subsection{Regularized Optimization and Fenchel-Young Loss}
    
            Since we consider a linear FOP in \Cref{eq:FOP}, the optimal solution $x^{\star}(\theta; u) \in \mathcal{X}^{\star}(\theta; u)$ may not be unique for a given $\theta$ and $u$ unless \Cref{asmp:almost-sure-uniqueness} holds. Moreover, $x^{\star}(\theta; u)$ is generally non-smooth with respect to $\theta$ for a fixed context $u \in \mathbb{R}^m$, implying that the sub-gradients $\partial_{\theta} x^{\star}(\theta; u)$ may not exist. Additionally, the decision error $D(\theta, \theta^{\star})$ in \Cref{eq:expected-decision-error} is non-convex in $\theta$ even if $x^{\star}(\theta; u)$ is convex in $\theta$ for every fixed context $u$, since $D(\theta, \theta^{\star})$ is not monotonic in $x^{\star}$ and convexity in $\theta$ is not preserved under composition \cite{boyd2004convex}. These challenges prevent the direct application of gradient-based methods to minimize $D(\theta, \theta^{\star})$ for $\theta$ over the parameter space $\Theta \subseteq \mathbb{R}^p$. 
    
            In view of these issues, we introduce a regularized version of the FOP in \Cref{eq:FOP} by adding a $1$-strongly-convex regularizer $\Omega$, for example, $\Omega(x) = \|x\|_2^2$, to the objective function:
            \begin{equation}\label{eq:regularized-FOP}
                \operatorname{R-FOP}(\theta^{\star}, u) := \max_{x \in \mathcal{X}(u)} \{h(\theta^{\star}; u)^{\top} x - \lambda \Omega(x)\}.
            \end{equation}
            Here $\lambda > 0$ is the regularization coefficient.
            The objective $\{h(\theta^{\star}; u)^{\top} x - \lambda \Omega(x)\}$ is strongly-concave in decision $x$, and the optimal solution is therefore unique:
            \begin{equation}\label{eq:regularized-decision}
                x_{\lambda}^{\star}(\theta^{\star}; u) := \argmax_{x \in \mathcal{X}(u)} \{h(\theta^{\star}; u)^{\top} x - \lambda \Omega(x)\}.
            \end{equation}
            We next show the regularized optimization is well-behaved. 
    
    
            \begin{proposition}[Properties of R-FOP]\label{prop:properties-regularized-FOP}
                Fix the context $u$. (1) \textbf{Unique Solution}: for any given parameter $\theta$ and regularization parameter $\lambda > 0$, the optimal solution of $\operatorname{R-FOP}(\theta, u)$ defined in \Cref{eq:regularized-decision} is unique. (2) \textbf{Lipschitz in parameter}: for fixed regularization parameter $\lambda > 0$, if $h(\theta; u)$ is $L$-Lipschitz in $\theta$ for any context $u$, then $x_{\lambda}^{\star}(\theta; u)$ is $\frac{L}{\lambda}$-Lipschitz with respect to $\theta$. (3) \textbf{Convergence}: for any $\theta \in \Theta$, we have $x_{\lambda}^{\star}(\theta; u) \rightarrow x^{\star}(\theta; u) $ as $\lambda \rightarrow 0$.
            \end{proposition}
    
            \Cref{prop:properties-regularized-FOP} suggests that the optimal solution of R-FOP is Lipschitz continuous in parameter $\theta$, and when $\lambda \rightarrow 0$, the optimal regularized solution $x_{\lambda}^{\star}(\theta; u)$ also converges to the optimal unregularized solution $x^{\star}(\theta; u)$ of FOP in \Cref{eq:FOP}. Accordingly, when $\lambda$ is small, the regularized FOP provides a good approximation for the original FOP.

    
            We now consider the suboptimality loss for the R-FOP with regularization parameter $\lambda > 0$, which maintains the convexity and smoothness of regularized objective in \Cref{eq:regularized-FOP} with respect to $\theta$. Define $V_{\lambda}(\theta; u, x) := h(\theta; u)^{\top} x - \lambda \Omega(x)$ as the value function of the regularized optimization in \Cref{eq:regularized-FOP}. The suboptimality loss associated with $\operatorname{R-FOP}(\theta, u)$ is
            \begin{equation}\label{eq:Fenchel-Young-loss}
                L_{\lambda}(\theta; u, y) := \max_{x \in \mathcal{X}(u)} V_{\lambda}(\theta; u, x) - V_{\lambda}(\theta; u, y).
            \end{equation}
            When $\lambda = 0$, $L_{\lambda}(\theta; u, y)$ reduces to the suboptimality loss in \Cref{eq:Sub-loss}. 
            
            Notably, the $L_{\lambda}(\theta; u, y)$ is precisely the Fenchel-Young (FY) loss defined on the region $\mathcal{X}(u)$ with observed decision $y$ and the regularizer $\lambda \Omega$ \citep{blondel2020learning,blondel2022learning}. 
            The FY loss was originally introduced for structured prediction (SP) tasks, which often involved implementing $\argmax$ or $\operatorname{rank}$ on a score function over a simplex constraint to generate the prediction. The inverse optimization therefore generlizes the SP tasks, as $\argmax$ and $\operatorname{rank}$ over simplex constraints therein can be formulated as linear optimization with convex constraints. Thus, our work extends the FY loss beyond SP to IOP.

\vskip 0.2in
    
        \subsection{Statistical Property of Fenchel-Young Loss}
    
            In \Cref{sec:formulation}, we show that the excess distance risk $R_{\rm DIST}(\theta) - \inf_{\theta} R_{\rm DIST}(\theta)$ is indeed equivalent to the decision error $D(\theta, \theta^{\star})$.
            By adapting the result from \citet[Proposition 3]{blondel2022learning}, we now establish the connection between the excess Fenchel-Young risk $R_{\lambda}(\theta) - \inf_{\theta} R_{\lambda}(\theta)$ and $D(\theta, \theta^{\star})$ in this part.
    
            \begin{theorem}[Calibration Bound]\label{thm:calibration-bound}
                Suppose \Cref{asmp:noisy-decisions,asmp:almost-sure-uniqueness} hold, and $\Omega$ is $1$-strongly convex. We have, for any $\lambda > 0$ and $\theta \in \mathbb{R}^p$, 
                \begin{equation*}
                    D(\theta, \theta^{\star}) = R(\theta) - R(\theta^{\star}) \leq 2 \mathbb{E}\big[\|x_{\lambda}^{\star}(\theta; U) - x^{\star}(\theta; U)\|_2^2 \big] + \frac{4}{\lambda} \big[ R_{\lambda}(\theta) - \inf_{\theta} R_{\lambda}(\theta) \big].
                \end{equation*} 
            \end{theorem}
            The term $\mathbb{E}[\|x_{\lambda}^{\star}(\theta; U) - x^{\star}(\theta; U)\|_2^2]$, represents the \emph{regularization error} induced by coefficient $\lambda$ for fixed $\theta$, and it vanishes as $\lambda \to 0$ according to \Cref{prop:properties-regularized-FOP} property (3). We discuss the control of regularization error and the excess FY risk in \Cref{sec:theoretical-analysis}.
    
            The key insight of \Cref{thm:calibration-bound} is that, when $\lambda$ is small, the excess FY risk $R_{\lambda}(\theta; u) - \inf_{\theta} R_{\lambda}(\theta; u)$ approximately provides an upper bound for the excess distance risk $R(\theta) - R(\theta^{\star})$. Therefore, the Fenchel-Young loss $L_{\lambda}(\theta; u, y)$ can be viewed as a convex and smooth \textbf{approximate surrogate loss} for the distance loss $L_{\rm DIST}(\theta; u, y)$ introduced by \citet{aswani2018inverse}. Therefore, minimizing the Fenchel-Young risk
            $$
                \min_{\theta \in \Theta} ~ R_{\lambda}(\theta) := \mathbb{E}[L_{\lambda}(\theta; U, Y)]
            $$
            leads to approximate minimization of the decision error $D(\theta, \theta^{\star})$, and therefore may provide a reasonable approximation of the unknown true parameter $\theta^{\star}$ of FOP in \Cref{eq:FOP}.        
    

\vskip 0.2in
            
        \subsection{Empirical Risk Minimization with Fenchel-Young Loss}
    
            In this part, we concretely discuss the properties FY loss following \citet{blondel2020learning,blondel2022learning}, which lay the foundation for designing a stochastic gradient descent (SGD) algorithm for minimizing the empirical Fenchel-Young risk:
            \begin{equation}\label{eq:IOP-FY}
                \hat{\theta}_{\lambda} \leftarrow \min_{\theta \in \Theta} ~ \widehat{R}_{\lambda}(\theta) := \frac{1}{n} \sum_{i=1}^{n} L_{\lambda}(\theta; U_i, Y_i).
            \end{equation}
            \citet{blondel2020learning,blondel2022learning} show that the Fenchel-Young loss $L_{\lambda}(\theta; u, y)$ is convex with respect to the parameter $\theta$ and admits a gradient under \Cref{asmp:linear-cost-function}, enabling the use of SGD.
    
            \begin{proposition}[\citet{blondel2022learning}, Proposition 3]\label{prop:properties-Fenchel-Young-loss}
                Suppose \Cref{asmp:linear-cost-function} holds.
                Fix $u \in \mathbb{R}^m$, $y \in \mathbb{R}^d$ and $\lambda > 0$. 
                (1) \textbf{Continuity}: $L_{\lambda}(\theta; u, y)$ is continuous in $\theta$ for any given $u$ and $y$, if the cost function $h(\theta; u)$ is continuous in $\theta$ for each given $u$.
                (2) \textbf{Convexity}: $L_{\lambda}(\theta; u, y)$ is convex in $\theta$ for any given $u$ and $y$. (3) \textbf{Gradient}: 
                By Danskin's theorem, we have $\nabla_{\theta} L_{\lambda}(\theta; u, y) = (\frac{\partial h(\theta; u)}{\partial \theta})^{\top} \left( x_{\lambda}^{\star}(\theta; u) - y \right)$, where $\frac{\partial h(\theta; u)}{\partial \theta}$ is the Jacobian matrix of $h$.
            \end{proposition}
    
            \Cref{prop:properties-Fenchel-Young-loss} is crucial for applying the SGD for solving the empirical Fenchel-Young risk $\widehat{R}(\theta)$ in \Cref{eq:IOP-FY}. The convexity of Fenchel-Young loss $L_{\lambda}(\theta; u_i, Y_i)$ with respect to $\theta$ leads to the convexity of empirical Fenchel-Young risk $\widehat{R}(\theta)$, as the non-negative sample averaging of convex functions is still convex \citep{boyd2004convex}. 
            The gradient of the empirical Fenchel-Young risk can be efficiently computed as the batch average of the gradients of the Fenchel-Young loss in each data point $(Y_i, U_i)$.
            We show the empirical FY risk minimization in \Cref{alg:FY-IOP-SGD}. 

\begin{algorithm}[t]
        \caption{Inverse Optimization with FY Loss (SGD)}
        \label{alg:FY-IOP-SGD}
        \begin{algorithmic}[1]
            \Require Noisy data $\{(U_i, Y_i)\}_{i=1}^{n}$, learning rate $\eta > 0$, batch size $m$, maximum iteration $T > 0$, tolerance $\epsilon > 0$ and parameter space $\Theta \subseteq \mathbb{R}^{p}$.
            \State Initialize parameter $\theta^{(0)} \in \Theta$.
            \For{$t = 0, 1, \dots, T$}
                \State Sample a batch $\mathcal{B}_t$ with size $m$ from $\{(U_i, Y_i)\}_{i=1}^{n}$.
                \State Compute the optimal solution $x_{\lambda}^{\star}(\theta^{(t)}; U_i), i \in \mathcal{B}_t$.
                \State Compute the batch-gradient of the empirical FY risk:
                $$
                    \nabla_{\theta} \widehat{R}_{\lambda}^{(t)} = \frac{1}{|\mathcal{B}_t|} \sum_{i \in \mathcal{B}_t} \left[\frac{\partial h(\theta; U_i)}{\partial \theta}\right]^{\top} (x_{\lambda}^{\star}(\theta^{(t)}; U_i) - Y_i).
                $$
                \State Update the parameter using the gradient-descent:
                $$
                    \theta^{(t+1)} = \theta^{(t)} - \eta \nabla_{\theta} \widehat{R}_{\lambda}(\theta^{(t)}).
                $$
                \If{$\| \nabla_{\theta} \widehat{R}_{\lambda}(\theta^{(t)}) \| \leq \epsilon$}
                    \State \textbf{break}.
                \EndIf
            \EndFor
            \State Return the estimated parameter $\hat{\theta}_{\lambda} \leftarrow \theta^{(T)}$.
        \end{algorithmic}
    \end{algorithm}

\vskip 0.3in
    
    \section{Theoretical Analysis}\label{sec:theoretical-analysis}
    
        In this section, we analyze the expected decision error $D(\hat{\theta}_{\lambda}, \theta^{\star})$ induced by the estimated parameter $\hat{\theta}_{\lambda}$ from \Cref{alg:FY-IOP-SGD}.
        We then provide theoretical guarantees for the regret $\operatorname{Reg}(\hat{\theta}_{\lambda}, \theta^{\star})$ and the parameter error $E(\hat{\theta}_{\lambda}, \theta) := \|\hat{\theta}_{\lambda} - \theta^{\star} \|_1$. 
        Detailed proofs are provided in supplemental materials.
        Before we proceed to the analysis, we first make the following assumption for the parameter space $\Theta$ in \Cref{alg:FY-IOP-SGD}. 
    
        \begin{assumption}[Restricted Parameter Space]\label{asmp:restricted-parameter-space}
            Suppose the parameter space $\Theta \subseteq \mathbb{R}^d$ satisfies that, (1) \textbf{Lower Bounded Norm}: for fixed context $u \in \mathbb{R}^m$, the norm of cost function $h(\theta; u)$ is lower bounded by a positive constant, say $b(u) > 0$, i.e., $\|h(\theta; u)\|_2 \geq b(u)$ for all $\theta \in \Theta$. (2) \textbf{Well-Specification}: for fixed context $u \in \mathbb{R}^m$, the parameter space $\Theta$ contains a parameter $\theta_{\lambda}^{\star}(u)$ such that $x_{\lambda}^{\star}(\theta_{\lambda}^{\star}(u); u) = \mathbb{E}[Y \mid u]$.
        \end{assumption}
    
         \Cref{asmp:restricted-parameter-space} could be a strong assumption. We provide some examples to illustrate the assumption in the supplemental materials.

        \subsection{Decision Error}\label{sub:decision-error}
    
            In this section, we provide a modified calibration bound for the expected decision error $D(\hat{\theta}_{\lambda}, \theta)$ under \Cref{asmp:restricted-parameter-space}, extending \Cref{thm:calibration-bound}. The key difference is that, the excess FY risk $R_{\lambda}(\theta) - \inf_{\theta \in \Theta} R_{\lambda}(\theta)$ is now evaluated on the restricted parameter space $\Theta$. 
    
            \begin{theorem}[Restricted Calibration Bound]\label{thm:restricted-parameter-space}
                Suppose \Cref{asmp:noisy-decisions,asmp:almost-sure-uniqueness,asmp:restricted-parameter-space} hold. The function $\Omega$ is $1$-strongly convex. For any $\lambda > 0$ and $\theta \in \Theta$, we have
                \begin{equation*}
                    D(\theta, \theta^{\star}) = R(\theta) - R(\theta^{\star}) \leq 2 \underbrace{ \mathbb{E}\big[\|x_{\lambda}^{\star}(\theta; U) - x^{\star}(\theta; U)\|_2^2 \big] }_{\text{Regularization Error}} + \frac{4}{\lambda} \underbrace{ \big[ R_{\lambda}(\theta) - \inf_{\theta \in \Theta} R_{\lambda}(\theta) \big] }_{\text{Excess FY Risk}}.
                \end{equation*} 
            \end{theorem}
            
\vskip 0.1in
            
            \noindent \textbf{Regularization Error}~~ For each fixed $\theta \in \Theta$, the regularization error measures the average distance of the regularized decision $x_{\lambda}^{\star}(\theta; U)$ and the unregularized decision $x^{\star}(\theta; U)$. 
            We have shown in \Cref{prop:properties-regularized-FOP} that $\lambda \rightarrow 0$, the expected regularization error vanishes. This suffices for establishing the consistency in the expected decision error.

            If we hope to establish non-asymptotic error bound, we need to further characterize how the regularization error $\mathbb{E}[\|x_{\lambda}^{\star}(\theta; U) - x^{\star}(\theta; U)\|_2^2]$ scales with the regularization coefficient $\lambda$.
            This requires analyzing the specific structure of the FOP, especially its feasible region $\mathcal{X}(u)$. In the following, we restrict our discussion to a special case, where the FOP is constrained in a ball that contains the origin, and the regularization function is $\Omega(x) = 1/2\|x\|_2^2$. We leave the general case for future work.
    
            \begin{theorem}[Lipschitzness of Regularized Solution under Ball Constraint]\label{thm:Lipschitz-regularized-oracle}
                Fix context $u \in \mathbb{R}^m$ and parameter $\theta \in \mathbb{R}^p$. Suppose the feasible set $\mathcal{X}(u) := \{x \in \mathbb{R}^d \mid \|x\|_2 \leq a(u)\}$ is a ball constraint with radius $a(u) > 0$, and the regularization function is $\Omega(x) = 1/2\|x\|_2^2$. If the critic value $\lambda_c(\theta; u) := \|h(\theta; u)\| / a(u) $ is positive, we then have 
                $$
                \left\{
                \begin{aligned}
                    & \|x_{\lambda}^{\star}(\theta; u) - x^{\star}(\theta; u)\|_2 = 0, ~ 0 \leq \lambda \leq \lambda_c(\theta; u), \\
                    & \|x_{\lambda}^{\star}(\theta; u) - x^{\star}(\theta; u)\|_2 \leq \lambda \frac{a(u)^2 }{2 \|h(\theta; u)\|_2}, ~ \lambda > \lambda_c(\theta; u).
                \end{aligned}
                \right.
                $$
            \end{theorem}
            \Cref{thm:Lipschitz-regularized-oracle} is essential to establish the bound of the regularization error for a fixed parameter $\theta \in \Theta$. However, the empirical minimizer $\hat{\theta}_{\lambda} \in \Theta$ depends on data and is therefore random. By assuming \Cref{asmp:restricted-parameter-space}, we can establish the uniform bound over $\Theta$:
            \begin{equation}\label{eq:regularization-error-bound}
                \mathbb{E}[\|x_{\lambda}^{\star}(\hat{\theta}_{\lambda}; U) - x^{\star}(\hat{\theta}_{\lambda}; U) \|_2^2] \leq \lambda^2 \mathbb{E}[a(U)^4 / (4 b(U)^2) ] .
            \end{equation}
            Therefore, as $\lambda \rightarrow 0$, the regularization error goes to zero, which confirms the result in \Cref{prop:properties-regularized-FOP}.
    
\vskip 0.1in
    
            \noindent \textbf{Excess Fenchel-Young Risk} ~~ To provide the generalization error bound for the excess Fenchel-Young risk $R_{\lambda}(\theta) - \inf_{\theta \in \Theta} R_{\lambda}(\theta)$, we introduce the empirical Rademacher complexity $\widehat{\mathcal{R}}_{\lambda}(\Theta)$ as
            $$
                \widehat{\mathcal{R}}_{\lambda}(\Theta) := \mathbb{E}_{\sigma}\Big[ \sup_{\theta \in \Theta} \frac{1}{n} \sum_{i=1}^{n} \sigma_i \tilde{L}_{\lambda}(\theta; U_i, Y_i) \Big].
            $$
            Here $\tilde{L}(\theta; U_i, Y_i) := L_{\lambda}(\theta; U_i, Y_i) - \mathbb{E}[L_{\lambda}(\theta; U_i, Y_i)]$ is a centered Fenchel-Young loss, and $\sigma = \{\sigma_1, \dots, \sigma_n\}$ is a Rademacher sequence. The empirical Rademacher complexity $\widehat{\mathcal{R}}_{\lambda}(\Theta)$ measures the complexity of the parameter space $\Theta$ with respect to the empirical FY loss. By applying the empirical process theory, we can derive the following generalization bound for the excess FY risk.
    
            \begin{lemma}[Symmetrization Bound, \citet{wainwright2019high} - Theorem 4.10]\label{lm:Rademacher-complexity-bound}
                Suppose that the parameter space $\Theta$ is bounded. Then for every finite $u$ and $y$, the FY loss $L_{\lambda}(\theta; u, y)$ is also uniformly bounded over $\Theta$ by some constant (denoted by $c$). For any $\lambda > 0$, we have, with probability at least $1 - \delta$,
                $R_{\lambda}(\hat{\theta}_{\lambda}) - \inf_{\theta} R_{\lambda}(\theta) \leq 2 \widehat{\mathcal{R}}_{\lambda}(\Theta) + \sqrt{ \frac{2 c^2 \log(1/\delta)}{n} }$.
            \end{lemma}

            In the statistical learning theory  literature, there are standard approaches to bound the empirical Rademacher complexity for many common function classes. The empirical Rademacher complexity often goes to zero as the sample size $n$ increases. We refer to the details in supplemental materials 
            and \citet{shalev2014understanding,mohri2018foundations} for more details. 

\vskip 0.1in
            
            \noindent \textbf{Decision Error Bound}~~ Suppose \Cref{asmp:noisy-decisions,asmp:linear-cost-function,asmp:restricted-parameter-space} hold. If the feasible region $\mathcal{X}(u)$ in \Cref{eq:regularized-FOP} is a ball-constrained region,  we can establish the finite-sample bound for the expected decision error $D(\hat{\theta}_{\lambda}, \theta^{\star})$. Combining \Cref{thm:restricted-parameter-space,thm:Lipschitz-regularized-oracle,lm:Rademacher-complexity-bound}, we have
            \begin{equation*}
                D(\hat{\theta}_{\lambda}, \theta^{\star}) \leq \frac{\lambda^2 \mathbb{E}[a(U)^4]}{4 b^2} + \frac{4}{\lambda} \Big[ \widehat{\mathcal{R}}_{\lambda}(\Theta) + \sqrt{\frac{c^2 \log(1/ \delta)}{n}} \Big].
            \end{equation*}
            We may further optimize the right-hand side of inequality with respect to $\lambda > 0$ to get a tighter  bound on the decision error.

\vskip 0.2in
    
        \subsection{Decision Regret and Parameter Error}\label{sub:regret-parameter-error}
    
            In this part, we show that the theoretical guarantees for the decision error $D(\hat{\theta}, \theta^{\star})$ directly lead to guarantees for the decision regret and the parameter error.

\vskip 0.1in

            \noindent \textbf{Decision Regret}~~ With the Cauchy-Schwarz inequality, we establish the following bound for the decision regret $\operatorname{Reg}(\hat{\theta}_{\lambda}, \theta^{\star})$ in terms of the decision error $D(\hat{\theta}_{\lambda}, \theta^{\star})$.
    
            \begin{theorem}[Decision Regret Bound]\label{thm:decisioin-regret-bound}
                If the expected squared norm of the true cost function $h(\theta^{\star}; U)$ is bounded by some constant, say $\mathbb{E}[\|h(\theta^{\star}; U)\|_2^2] \leq B$, we have $\operatorname{Reg}(\hat{\theta}_{\lambda}, \theta^{\star}) \leq B D(\hat{\theta}_{\lambda}, \theta^{\star})$.
            \end{theorem}
    
            \Cref{thm:decisioin-regret-bound} indicates that the convergence of decision error $D(\hat{\theta}_{\lambda}, \theta^{\star})$  implies the convergence of decision regret $\operatorname{Reg}(\hat{\theta}_{\lambda}, \theta^{\star})$, and any bound for the decision error implies the corresponding bound for the decision regret.

\vskip 0.1in
            
            \noindent \textbf{Parameter Error} ~~ By assuming the identifiability condition proposed in \citet[IC Conditions]{aswani2018inverse}, we analyze the convergence of parameter error $E(\hat{\theta}_{\lambda}, \theta^{\star})$ in terms of the convergence of decision error $D(\hat{\theta}_{\lambda}, \theta^{\star})$.
    
            \begin{theorem}[Uniform Convergence of Parameter Error]\label{thm:uniform-convergence-parameter-error}
                Suppose \Cref{asmp:almost-sure-uniqueness} holds, and for all $\theta \in \Theta \backslash \theta^{\star}$, $\operatorname{dist}(\mathcal{X}^{\star}(\theta; u), \mathcal{X}^{\star}(\theta^{\star}; u)) > 0$. Then for any $\theta$ such that $x^{\star}(\theta; u) \rightarrow_p x^{\star}(\theta^{\star}; u)$, we have $\theta \rightarrow_p \theta^{\star}$.
            \end{theorem}

            \Cref{thm:uniform-convergence-parameter-error} establishes the statistical convergence of decision error $D(\hat{\theta}_{\lambda}, \theta^{\star})$ implies the statistical convergence of parameter error $E(\hat{\theta}_{\lambda}, \theta^{\star})$ to zero.

\begin{table}[t]
    \caption{Overview of the Synthetic Data Experiments.}
    \centering
    \vskip 0.1in
    \begin{tabular}{c|clll}
    \toprule
    Experiments & \multicolumn{4}{c}{Forward Optimization Problem} \\ 
    \midrule
    Example A & \multicolumn{4}{l}{$\min \{ (\theta + u)^{\top} x \mid x \succeq 0, \parallel x \parallel_1 \leq a \}$} \\
    Example B & \multicolumn{4}{l}{$\min \{ (\theta \circ u)^{\top} x \mid x \in [-1, 1]^{p} \}$} \\
    Example C & \multicolumn{4}{l}{$\min \{ (\theta + u)^{\top} x \mid x \in [-1, 1]^{p} \}$} \\
    Example D & \multicolumn{4}{l}{$\min \{ x^{\top}x - (\theta + u)^{\top} x \mid x \in [0, 1]^{p} \}$} \\
    Example E & \multicolumn{4}{l}{$\min \{ -(\theta + u)^{\top} x \mid \parallel x \parallel_2^2 \leq a^2 \}$} \\ \bottomrule
    \end{tabular}
    \label{tab:Numerical-experiments-overview}
    \end{table}

    \begin{table}[t]
    \caption{The experimental results for various methods under different metrics for \textbf{Example A} in the {Noisy Decision} setting, where the FOP is  defined as $\min \{ (\theta + u)^{\top} x \mid x \succeq 0, \|x\|_1 \leq a \}$. For the Fenchel-Young (FY) loss, we set $\lambda = 0.1$ and $\Omega(x) = 1/2\|x\|_2^2$.}
    \centering
    \vskip 0.1in
    \begin{tabular}[width=\linewidth]{c|cccc|cccc|cccc}
    \hline
    \multirow{2}{*}{Sample size} & \multicolumn{4}{c|}{Parameter Error} & 
    \multicolumn{4}{c|}{Decision Error} & \multicolumn{4}{c}{Regret} \\ 
    \cline{2-13} 
     & FY & SPA & KKA & VIA & FY & SPA & KKA & VIA & FY & SPA & KKA & VIA \\ 
     \hline
    50 & \textbf{1.30} & 2.46 & 4.25 & 3.81 & \textbf{3.03} & 4.53 & 7.66 & 3.90 & \textbf{0.05} & 0.20 & 0.77 & 0.08 \\
    100 & \textbf{1.17} & 1.35 & 4.58 & 4.12 & \textbf{2.47} & 2.94 & 7.79 & 3.17 & \textbf{0.03} & 0.06 & 0.84 & 0.05 \\
    300 & \textbf{0.53} & 2.10 & 4.87 & 4.42 & \textbf{2.05} & 2.66 & 8.22 & 2.19 & \textbf{0.02} & 0.08 & 0.93 & 0.03 \\
    500 & \textbf{0.53} & 2.44 & 4.92 & 4.52 & \textbf{1.64} & 2.69 & 8.32 & 1.88 & \textbf{0.01} & 0.10 & 0.95 & 0.03 \\
    1000 & \textbf{0.38} & 2.59 & 4.96 & 4.56 & \textbf{1.33} & 2.64 & 8.40 & 1.60 & \textbf{0.01} & 0.11 & 0.96 & 0.03 \\ \hline
    \end{tabular}
    \label{tab:example-A}
    \end{table}    
    
\begin{table}[t]
    \caption{The experimental results for various methods under different metrics for \textbf{Example B} in the {Noisy Decision} setting, where the FOP is  defined as $ \min \{ (\theta \circ u)^{\top} x \mid x \in [-1, 1]^{p} \}= \min \left\{\sum_{k=1}^{p} \theta_k  u_k x_k \mid x \in [-1, 1]^{p} \right\}$. We set $\lambda=0.1$ and $\Omega(x) = 1/2 \|x\|_2^2$.}
    \centering
    \vskip 0.1in
    \begin{tabular}{c|cccc|cccc|cccc}
    \hline
    \multirow{2}{*}{Sample size} & \multicolumn{4}{c|}{Parameter Error} & \multicolumn{4}{c|}{Decision Error} & \multicolumn{4}{c}{Regret} \\ \cline{2-13} 
     & FY & SPA & KKA & VIA & FY & SPA & KKA & VIA & FY & SPA & KKA & VIA \\ \hline
    50 & \textbf{2.76} & 5.50 & 5.50 & 5.50 & \textbf{0.00} & 19.99 & 19.99 & 0.00 & \textbf{0.00} & 2.75 & 2.75 & 0.00 \\
    100 & \textbf{2.42} & 5.50 & 5.50 & 5.50 & \textbf{0.00} & 19.99 & 19.99 & 0.00 & \textbf{0.00} & 2.76 & 2.76 & 0.00 \\
    300 & \textbf{2.25} & 5.50 & 5.50 & 5.50 & \textbf{0.00} & 20.00 & 20.00 & 0.00 & \textbf{0.00} & 2.75 & 2.75 & 0.00 \\
    500 & \textbf{2.27} & 5.50 & 5.50 & 5.50 & \textbf{0.00} & 20.01 & 20.01 & 0.00 & \textbf{0.00} & 2.76 & 2.76 & 0.00 \\
    1000 & \textbf{2.19} & 5.50 & 5.50 & 5.50 & \textbf{0.00} & 19.99 & 19.99 & 0.00 & \textbf{0.00} & 2.75 & 2.75 & 0.00 \\ \hline
    \end{tabular}
    \label{tab:Example-B}
\end{table}

\vskip 0.3in
    
\section{Experiments}\label{sec:experiments}
    
    We conducted a series of synthetic and real data experiments, and we briefly present some of the results in this section, with more details deferred to \Cref{sec: Experiment-details-appendix}.
     
\vskip 0.2in
    
    \subsection{Synthetic Data Experiments}
    
    We design a total of five synthetic data experiments, which encompass a variety of scenarios including box constraints, half-space constraints, ball constraints, linear objective functions and nonlinear objective functions, as shown in \Cref{tab:Numerical-experiments-overview}.  For evaluation, we report the computational time, parameter error, decision error, and regret. All the results are averaged over 100 experiment replications.
    
    Considering different forms of noisy observations, we examine three distinct settings for each example: \textbf{Noisy Decision} in \Cref{asmp:noisy-decisions}, \textbf{Noisy Objective} in \Cref{rmk:noisy-objective}, and a \textbf{Noiseless} setting (the true optimal decision of FOP is observed). For the former two settings, standard normal noises are added to the decision observations and the optimization objectives respectively. 
    
    In this section, we present the results for Examples A and B in the Noisy Decision setting. Specifically, we compare our method with other benchmark methods discussed in \Cref{sub:LossFunction}. These methods include the SPA (Semi-parametric methods), KKA (KKT-Approximated Loss), and VIA (Variational Inequality Loss) discussed in \Cref{sub:LossFunction}. The results are shown in \Cref{tab:example-A} and \Cref{tab:Example-B} respectively. Notably, both Example A and Example B involve linear optimization, for which the VIA loss suboptimality Loss is equivalent to suboptimality loss as discussed in \Cref{sub:LossFunction}. 
    More experimental results for other settings are available in \Cref{sub:synthetic-experiment}.

\vskip 0.1in
    
    \noindent \textbf{Results for Example A} ~~
    Our approach, leveraging the Fenchel-Young loss, demonstrates superior performance across key metrics, including parameter error, decision error, and regret. Notably, it achieves the lowest values for these metrics compared to other traditional inverse optimization methods. While most methods, with the exception of KKA, exhibit gradual convergence in decision error and regret as the sample size increases, their performance in parameter recovery remains consistently suboptimal. Even worse, the KKA method fails to converge in all of these metrics.


\vskip 0.1in

    \noindent \textbf{Results for Example B} ~~
    The FY loss and VIA demonstrate the best performance, achieving zero decision error and regret even with small sample sizes. In contrast, the KKA and SPA methods perform poorly with non-vanishing errors. Notably, the parameters $\theta$ estimated by SPA, KKA, and VIA often converge toward zero, causing problem degeneracy. Among these methods, VIA happens to  successfully recovere the correct sign of each element in $\theta$, hence achieving zero decision error and regret. However, their nearly zero estimates are not very meaningful. 
    
    As explained in the \Cref{sec:example-B-details}, 
    the degeneracy arises because, in the optimization models of KKA, VIA, and SPA, $\theta = \mathbf{0}$ coincidentally represents an optimal solution. 
    To alleviate the degeneracy issue, we have also tried introducing various regularization constraints, such as constraining the sum of $\theta$’s elements or the first element of $\theta$  to a nonzero value. However, these modifications did not lead to significant improvements so we omit them for brevity.


\vskip 0.1in
    
    \noindent \textbf{Computational Time} ~~ In most experiments, our method achieves the shortest computational time across all experiments with an order-of-magnitude improvement. More detailed results are available in \Cref{sec:computational-time-appendix}.

\vskip 0.1in
    
    \noindent \textbf{Regularization Parameter} ~~
    Our experiments demonstrate that a small value of $\lambda > 0$ is sufficient to achieve robust performance. As the sample size increases, reducing $\lambda$ further tends to improve results. For sufficiently large sample sizes, cross-validation can be used to select the optimal regularization parameter. A detailed analysis of how the regularization parameter $\lambda$ impacts the performance of the FY loss $L_{\lambda}$ is provided in \Cref{sec:sensitivity-analysis-appendix}.
    

\vskip 0.2in
    
    \subsection{Real data Experiments}

\begin{table}[t]
    \caption{
    The \textbf{Relative Regret Ratio} column reports percentages (e.g., $1.01$ represents a regret ratio of $1.01\%$). The \textbf{Period} column indicates the dataset’s time span in years (e.g., $1.5$ corresponds to one and a half years). The \textbf{Computational Time} column lists average running times in seconds (e.g., $6$ means an average of $6$ seconds). For the Fenchel-Young loss, we set $\lambda = 0.1$ and $\Omega(x) = 1/2\|x\|_2^2$.
    }
    \vskip 0.1in
    \centering
    \resizebox{\textwidth}{12mm}{
    \begin{tabular}{c|clccc|ccccc|ccccc}
    \hline
    \multirow{2}{*}{\begin{tabular}[c]{@{}c@{}}Period\\ (years)\end{tabular}} & \multicolumn{5}{c|}{Decision Error} & \multicolumn{5}{c|}{Relative Regret Ratio (\%)} & \multicolumn{5}{c}{Computing time (seconds)} \\ \cline{2-16} 
     & FY & \multicolumn{1}{c}{SPA} & KKA & VIA & MOM & FY & SPA & KKA & VIA & MOM & FY & SP & KKA & VIA & MOM \\ \hline
    0.5 & \textbf{2.06} & 6.80 & 5.41 & 10.72 & 3.58 & \textbf{1.01} & 29.51 & 8.30 & 34.51 & 4.87 & \textbf{6} & 120 & 66 & 59 & 438 \\
    1 & \textbf{2.00} & 6.72 & 6.61 & 11.77 & 4.09 & \textbf{0.97} & 29.55 & 12.39 & 41.51 & 4.96 & \textbf{8} & 287 & 152 & 107 & 647 \\
    1.5 & \textbf{1.97} & 6.77 & 7.09 & 7.33 & 4.66 & \textbf{0.94} & 29.50 & 12.81 & 30.95 & 4.21 & \textbf{9} & 513 & 242 & 185 & 901 \\
    2 & \textbf{1.89} & 6.78 & 7.07 & 7.17 & 4.64 & \textbf{0.84} & 29.28 & 12.40 & 30.41 & 3.26 & \textbf{11} & 839 & 281 & 256 & 1218 \\ \hline
    \end{tabular}
    }
    \label{tab:real-data-experiment-result}
    \end{table}
    
    We conducted experiments using a real-world dataset from Uber Movement (\url{https://movement.uber.com}), focusing on 45 census tracts in downtown Los Angeles. This dataset includes average travel times between neighboring tracts across 93 edges for 2018 and 2019. To construct an inverse optimization problem, we solve for a shortest path for each travel time observation in the dataset, which serves as the observed optimal solution $Y$. Accordingly, we consider a shortest path problem as the FOP. 
    
    In our experiment setting, the context vector $u \in \mathbb{R}^{12}$ includes features such as wind speed and visibility. The travel times vector for 93 edges is modeled as $\boldsymbol{\theta} u$, where $\boldsymbol{\theta} \in \mathbb{R}^{93 \times 12}$ is the unknown parameter matrix. The shortest path FOP is given by $\min_{x} \left\{ \sum_{k=1}^{93} (\boldsymbol{\theta}_k u) x_k \mid Ax = b \right\}$, where $\boldsymbol{\theta}_k $ is the k-th row of the parameter matrix $\boldsymbol{\theta}$ and $Ax = b$ is the flow constraints with $A \in \mathbb{R}^{45 \times 93}$ and $b \in \mathbb{R}^{45}$. The IOP task is to estimate $\boldsymbol{\theta}$ from observations of shortest path solutions. 
    For more experiment setup details, please see \Cref{sub:real-data-experiment}.
    
    In this real-data experiment, we continue to compare the FY loss approach with the benchmarks described in \Cref{sub:LossFunction}. 
    Additionally, we note that the Maximum Optimality Margin (MOM) method proposed by \citet{sun2023maximum} can also be applied in this setting, we also include it as a benchmark method for comparison. Readers interested in this method can refer to \Cref{sub:LossFunction}.

    We divide the dataset into four time spans: half a year (920 points), one year (1825 points), one and a half years (2735 points), and two years (3640 points). For each case, $60\%$ of the data is allocated for training and $40\%$ for testing. This random split is repeated across 30 independent experimental trials. The metrics collected include computational time, decision error, and relative regret ratio\footnote{The relative regret ratio quantifies the testing regret of each method relative to the clairvoyant shortest paths on the testing data, normalized by the average traveling time of the clairvoyant shortest paths.}, where the decision error and regret ratio are both calculated on the testing data. We summarize the experimental results in \Cref{tab:real-data-experiment-result} as below.

\vskip 0.1in

    \noindent \textbf{Decision Error}~~ The FY loss method achieves the lowest decision error across all time spans, significantly outperforming baseline approaches. This demonstrates its consistent ability to identify paths closely aligned with the true optimal solution. The MOM method ranks second, delivering near-optimal paths but falling short of our method’s accuracy. In contrast, SPA, VIA, and KKA exhibit substantially larger deviations from the optimal path.

\vskip 0.1in

    \noindent \textbf{Regret Ratio}~~  Consistent with its decision error performance, our method again surpasses all competitors, producing path travel times that tightly converge to the theoretical optimum. MOM remains the best competitor, maintaining a relatively low regret ratio, albeit still significantly higher than our approach. In contrast, KKA, SPA and VIA perform poorly, achieving even larger regret ratios. 

\vskip 0.1in

    \noindent \textbf{Computational Time}~~ Our method delivers superior computational efficiency, completing tasks in seconds compared to the minutes or hours required by other methods. Competing approaches (SPA, VIA, KKA, MOM) rely on solving data-intensive constrained optimization problems, where runtime scales fast with input size. In contrast, our framework leverages a finite sequence of gradient descent steps, paired with lightweight gradient computations, to achieve rapid convergence. This orders-of-magnitude speed advantage makes our method uniquely viable for data-driven IOP with large datasets.
    
\vskip 0.2in    

\section{Conclusion}
    In this paper, we build a connection between inverse optimization and the Fenchel-Young (FY) loss originally designed for structured prediction.
    We propose a FY loss approach to data-driven inverse optimization with an efficient gradient-based optimization algorithm and rigorous theoretical guarantees. 
    We demonstrate that this approach significantly outperforms many existing data-driven inverse-optimization methods in extensive simulations and real-data experiments. 
    We hope this paper provides a new perspective on both FY loss and 
    inverse optimization, opening more future research on FY loss for inverse optimization. 

\clearpage

\bibliographystyle{informs2014} %
\bibliography{inverse-optimization} %

\ECSwitch
\EquationsNumberedBySection

\section{Supplementary for \Cref{sec:formulation}}\label{sec:supplementary-formulation}

    \subsection{Some Examples of Linear Forward Optimization}

        \begin{example}[Consumer Utility Models]\label{ex:consumer-demand-models}
            The consumer demand model are widely used in characterize the consumer utility in electricity market \citep{saez2016data,saez2017short,fernandez2021forecasting,fernandez2021inverse}. For each time period $t$, the model can be formulated as:
            \begin{equation}\label{eq:consumer-demand-model}
                \max_{y}\{u_t(\theta; y) - p_t y \mid l_t \leq y \leq r_t\},
            \end{equation}
            where $y$ is the electricity consumption level, restricted by $l_t$ and $r_t$, $p_t$ is the price at time $t$, and $\theta$ is the unknown parameter that characterizes the consumer utility function $u_t(\theta; y)$. 
            As summarized in \citet{chan2023inverse}, the utility $u_t(\theta; y)$ can be chosen as a piecewise linear function in $\theta$:
            $$
                u_t(\theta; y) = \left\{
                \begin{aligned}
                    & \theta_{1,t} y, ~~ && \text{if} ~ y < s, \\
                    & \theta_{1,t} y + \cdots + \theta_{i, t} (y - (i - 1)s), ~~ && \text{if} ~ (i-1)s \leq y \leq i s, 2 \leq i \leq d.
                \end{aligned}
                \right.
            $$
            By letting $h(\theta_t; u) = \theta_t$, $\theta_t = (\theta_{1, t}, \dots, \theta_{d, t})^{\top} \in \mathbb{R}^d $, $\bar{y} = (y \mathbb{I}_{y < s}, \dots, y \mathbb{I}_{(d-1)s \leq y < s} )^{\top} \in \mathbb{R}^d $,  and $g(u, \bar{y}) = (l_t - y, y - r_t)^{\top} \in \mathbb{R}^2$, the FOP in \Cref{eq:FOP} exactly captures the consumer demand model in \Cref{eq:consumer-demand-model}.
        \end{example}

        \begin{example}[Vehicle Routing Problem]\label{ex:vehicle-routing-problem}
            \citet{chen2021inverse} consider a vehicle routing problem with limited capacity. The FOP is formulated as a mixed-integer linear optimization model, as summarized in \citet{chan2023inverse}:
            $$
                \min_{x, y}\{\theta^{\top} x \mid Ax + By \geq d, \pi \leq y \leq Q \mathbf{1}, x \in \{0, 1\}^p \}.
            $$
            Here, the variables $x$ and $y$ represent the arc choices and the amount of goods each vehicle carries, respectively. The unknown parameter $\theta$ is the cost vector of the arcs. The problem is constrained by the vehicle capacity $Q$ and flow balance constraints. In this example, the objective is bilinear in the parameter $\theta$ and decision $x$, and the constraint is convex, fitting our formulation in \Cref{eq:FOP}.
        \end{example}

\vskip 0.2in

    \subsection{Proof of \Cref{lm:excess-distance-risk}}

        \begin{proof}{Proof of \Cref{lm:excess-distance-risk}}
            By \Cref{asmp:almost-sure-uniqueness}, we have
            $$
                R(\theta) = \mathbb{E}[\min_{x \in \mathcal{X}^{\star}(\theta; u)} \|Y - x\|_2^2] = \mathbb{E}[\|Y - x^{\star}(\theta; U)\|_2^2].
            $$
            The excess risk is then
            $$
            \begin{aligned}
                R(\theta) - R(\theta^{\star}) &= \mathbb{E}\big[ \|Y - x^{\star}(\theta; U)\|_2^2 - \|Y - x^{\star}(\theta^{\star}; U)\|_2^2 \big] \\
                &= \mathbb{E}\big[ \|x^{\star}(\theta; U)\|_2^2 \big] - \mathbb{E}\big[ \|x^{\star}(\theta^{\star}; U)\|_2^2 \big] + 2 \mathbb{E}\big[ \mathbb{E}[Y \mid U]^{\top} (x^{\star}(\theta^{\star}; U) - x^{\star}(\theta; U) ) \big]
            \end{aligned}
            $$
            By \Cref{asmp:noisy-decisions}, we have $\mathbb{E}[Y \mid U] = x^{\star}(\theta^{\star}; U)$ almost surely. Therefore, the excess risk can be further simplified as
            $$
            \begin{aligned}
                R(\theta) - R(\theta^{\star}) &= \mathbb{E} \big[\|x^{\star}(\theta; U)\|_2^2 \big] - \mathbb{E}\big[ \|x^{\star}(\theta^{\star}; U)\|_2^2] + 2 \mathbb{E}[ \|x^{\star}(\theta^{\star}; U)\|_2^2 \big] - 2 \mathbb{E}\big[ x^{\star}(\theta^{\star}; U)^{\top} x^{\star}(\theta; U) \big] \\
                &= \mathbb{E}\big[ \|x^{\star}(\theta; U) - x^{\star}(\theta^{\star}; U)\|_2^2 \big],
            \end{aligned}
            $$
            where the expectation is taken over the marginal distribution of $U$. This completes the proof.

        \end{proof}  

\vskip 0.2in

    \subsection{Loss Functions for Data-driven Inverse Optimization}\label{sub:loss-functions}

        \paragraph{KKT-Approximated Loss} \citet{keshavarz2011imputing} consider a noiseless decision setting and propose a loss function that measures the violation of the Karush-Kuhn-Tucker (KKT) conditions:  
        \begin{equation*}\label{eq:KKA-loss}
            L_{\rm KKA}(\theta; u, y) := \min_{\lambda \succeq 0} \{L_{\rm ST}(\theta; u, y, \lambda) + L_{\rm CS}(u, y, \lambda)\},
        \end{equation*}
        where $\lambda$ is the dual variable and the two loss functions corresponds to the \emph{stationary condition} and \emph{complementary slackness} as follow:
        $$
        \begin{aligned}
            & L_{\rm ST}(\theta; u, y, \lambda) := \big\| h(\theta; u) + \sum_{j=1}^{m} \lambda_j \nabla_x g_j(x, u) \big\|_2, \\
            & L_{\rm CS}(u, y, \lambda) := \big\| \big(\lambda_1 g_1(x, u), \dots, \lambda_m g_m(x, u) \big) \big\|_2.
        \end{aligned}
        $$
        In our setting with noisy observations, minimizing the IOP risk with $L_{\rm KKA}$ is challenging, as the decision may be entirely infeasible and thus significantly violate the KKT conditions. The experiments in \Cref{sec:experiments} show that the KKT loss is sensitive to noise and may not converge to a reasonable estimation even when the sample size is large.

        \paragraph{Suboptimality Loss} 

        \citet{chan2014generalized,chan2018trade} study the inverse optimization of multi-objective linear programming with single observation by minimizing the absolute duality gap between the primal objective value induced by observed decision and the dual objective value induced by the parameter.
        \citet{barmann2017emulating} consider the suboptimality loss in the online learning setting.
        \citet{chan2019inverse} extend this framework to single noisy observation and proposed a goodness-of-fit metric for IOP. 
        They proposed to estimate the parameters by minimizing the average duality gap between the objective function of observed decisions and the dual optimal objective $v(\theta; \lambda_i, U_i)$ induced by the parameter:
        \begin{equation*}
            \min_{\theta, \epsilon_i, \lambda_i \geq 0, i = 1, \dots, n} \frac{1}{n} \sum_{i=1}^{n} \epsilon_i ~ \text{s.t.} ~ v(\theta; \lambda_i, U_i) - h(\theta; U_i)^{\top} Y_i  \leq \epsilon_i. \\
        \end{equation*}
        \citet{chan2023inverse} points out \Cref{eq:duality-gap} is equivalent to minimizing the empirical risk $\widehat{R}(\theta)$ in \Cref{eq:empirical-risk-minimization} with the loss replaced by the suboptimality loss as follows: 
        \begin{equation*}
            L_{\rm SUB}(\theta; u, y) := \max_{x \in \mathcal{X}(u)} h(\theta; u)^{\top} x - h(\theta; u)^{\top} y.
        \end{equation*}
        Although the optimization problem in this approach is convex and tractable, this duality-based approach still faces several challenges. (1) If there exists a parameter $\theta \in \Theta$ such that $h(\theta; u) = 0$, the suboptimality loss $L_{\rm SUB}$ also reaches zero, causing the inverse optimization problem to degenerate. (2) As the sample size and decision dimension increase, the existing dual formulation in \Cref{eq:duality-gap} for solving empirical risk induced by suboptimality loss becomes computationally challenging. (3) It is unknown that whether the empirical risk minimization with the suboptimality loss leads to the theoretical control for the decision error we target in \Cref{eq:expected-decision-error}.

        \paragraph{Variational Inequality Loss} 
        Inspiring by the \emph{variational inequality}, \citet{bertsimas2015data} study the variational inequality approximation (VIA) loss function that measures the optimality margin for the FOP with a differentiable objective function $f(\theta; u, x)$:
        $$
            \nabla_x f(\theta; u, x)^{\top} (x - x^{\star}(\theta; u)) \geq 0,
        $$
        where $f(\theta; u, x)$ is the objective function of the FOP. For an observed decision $y \in \mathbb{R}^d$, the VIA loss is defined as:
        \begin{equation*}
            L_{\rm VIA}(\theta; u, y) = \max_{x \in \mathcal{X}(u)} \nabla_y f(\theta; u, y)^{\top} (x - y).
        \end{equation*}
        As noted by \citet{chan2023inverse}, for linear FOPs, the VIA loss is equivalent to the suboptimality loss in \Cref{eq:Sub-loss}, as:
        $$
            L_{\rm VIA}(\theta; u, y) = |\max_{x \in \mathcal{X}} h(\theta; u)^{\top} (x - y)| = L_{\rm SUB}(\theta; u, y).
        $$
        \citet{bertsimas2015data} propose solving the VIA loss via strong duality formulation, which facing the same computational challenges as the suboptimality loss. The VIA loss is also sensitive to noise, as observed decisions may fall outside the feasible region, leading to large VIA losses.

        \paragraph{Distance Loss}  

        \citet{aswani2018inverse} propose an IOP algorithm that directly minimizes the expected decision error $E(\theta, \theta^{\star})$ over $\Theta$ in \Cref{eq:expected-decision-error} by considering minimized the squared distance disk as shown in \Cref{eq:IOP-risk}.
        \begin{equation*}
            L_{\rm DIST}(\theta; u, y) := \min_{x \in \mathcal{X}^{\star}(\theta;u) } \big\|y - x \big\|_2^2.
        \end{equation*}
        By enforcing strong duality, \citet{aswani2018inverse} reformulate the empirical IOP risk with distance loss as a convex optimization problem:
        \begin{equation*}
            \begin{aligned}
                \widehat{R}_{\rm DIST}(\theta) := \min_{x_i, \lambda_i} ~~ & \frac{1}{n} \sum_{i=1}^{n} \|Y_i - x_i\|_2^2 \\
                \text{s.t.} ~~ & v(\theta; \lambda_i, U_i) - h(\theta; U_i)^{\top} x_i \leq 0, && \forall ~ i \in [n] \\
                & g(U_i, x_i) \preceq 0, && \forall ~ i \in [n] \\
                & \lambda_i \geq 0, && \forall ~ i \in [n].
            \end{aligned}
        \end{equation*}
        Here, $v(\theta; \lambda_i, U_i)$ is the Lagrangian dual function of the FOP, and $\lambda_i$ is the dual variable for each observed decision $Y_i$. 
        Note that, even when \Cref{asmp:almost-sure-uniqueness} holds, directly minimizing the expected decision error $D(\theta, \theta^{\star})$ or the IOP risk $R(\theta)$ in \Cref{eq:IOP-risk} over $\Theta$ presents significant challenges. The primary difficulty arises from the non-convexity of the objective function  $\|x^{\star}(\theta; U) - x^{\star}(\theta^{\star})\|_2^2$ , making it difficult to find the global minimum. In addition, $D(\theta, \theta^{\star})$ may not be continuous in $\theta$, which means the subgradient $\partial_{\theta}D(\theta, \theta^{\star})$ may not exist, preventing the use of gradient-based optimization methods. To address this, \citet{aswani2018inverse} propose an enumeration algorithm that evaluates  $R(\theta)$ over a discretized parameter space $\Theta$. However, this approach becomes computationally expensive in high-dimensional settings. Moreover, computing $R(\theta)$ involves solving a convex optimization problem, whose complexity scales significantly with the sample size and decision dimension, further exacerbating computational challenges.

        \paragraph{Semi-parametric Methods} 

        To overcome the limitations of squared distance loss with large problems, \citet{aswani2018inverse} proposed a semi-parametric (SPA) method that first uses the Nadaraya-Watson (NW) estimator to denoise the observed decisions and then projects them onto the feasible region, denoted as $\tilde{X}_i$. 
        \citet{aswani2018inverse} then propose to solve \Cref{eq:duality-gap} by replacing $Y_i$ with $\tilde{X}_i$ for each $i \in [n]$, that is,
        \begin{equation}
            \begin{aligned}
                \min_{\theta, \epsilon_i, \lambda_i} \frac{1}{n} \sum_{i=1}^{n} \epsilon_i ~~ \text{s.t.} ~ & v(\theta; \lambda_i, U_i) - h(\theta; U_i)^{\top} \tilde{X}_i \leq \epsilon_i, \\
                & \theta \in \Theta, ~~ \lambda_i \geq 0, ~~ \forall ~ i \in [n].
            \end{aligned}
        \end{equation}
        However, the NW estimator performs poorly with high-dimensional context space, and the dual formulation is also computationally expensive when the sample size and dimension grow. Meanwhile, such estimate-then-project approach separates the de-noising process from the structure of FOP, potentially resulting in suboptimal performance.
        

        \paragraph{Maximum Optimality Margin} \citet{sun2023maximum} considered an inverse optimization method for which the FOP is a linear program (LP). Specifically, they analyzed the following LP from the perspective of the simplex method, assuming that the observed solution $x^*=\left(x_1, \ldots, x_n\right)^{\top}$ is noise-free.
        \vskip 0.in
        $$
        \begin{aligned}
            \operatorname{LP}(c, A, b):= \min ~~ & c^{\top} x \\
            \text {s.t.} ~~ & A x=b, x \geq 0
        \end{aligned}
        $$
        where $c \in \mathbb{R}^n, A \in \mathbb{R}^{m \times n}$, and $b \in \mathbb{R}^m$. Based on the simplex method, they define the optimal basis $\mathcal{B}^*$ and its complement $\mathcal{N}^*$ as follows,
        $$
            \mathcal{B}^*:=\left\{i: x_i^*>0\right\}, \quad \mathcal{N}^*:=\left\{i: x_i^*=0\right\} .
        $$
        For a set $\mathcal{B} \subset[n]$, they use $A_{\mathcal{B}}$ to denote the submatrix of $A$ with column indices corresponding to $\mathcal{B}$ and $c_{\mathcal{B}}$ to denote the subvector with corresponding dimensions. Therefore, we have $\mathcal{B}=\mathcal{B}^*$ if and only if 
        \begin{equation}\label{eq:MOM}
            c_{\mathcal{N}}^{\top}-c_{\mathcal{B}}^{\top} A_{\mathcal{B}}^{-1} A_{\mathcal{N}} \geq 0
        \end{equation}
        Suppose there is a linear mapping from the covariates $z_t$ to the objective vector $c_t$, i.e., $g\left(z_t ; \Theta\right):=\Theta z_t$
        The Maximum Optimality Margin (MOM) method solves the following optimization problem.
        $$  
        \begin{aligned}
            \hat{\Theta} := \argmin_{\Theta \in \mathcal{K}} & ~~ \frac{\lambda}{2}\|\Theta\|_2^2+\frac{1}{T} \sum_{t=1}^T\left\|s_t\right\|_1 \\
            \text {s.t.} & ~~ \hat{c}_t=\Theta z_t, \quad t=1, \ldots, T, \\
            & \hat{c}_{t, \mathcal{N}_t^*}^{\top}-\hat{c}_{t, \mathcal{B}_t^*}^{\top} A_{t, \mathcal{B}_t^*}^{-1} A_{t, \mathcal{N}_t^*} \geq 1_{\left|\mathcal{N}_t^*\right|}-s_t, ~ t=1, \ldots, T.
        \end{aligned}
        $$
        where the decision variables are $\Theta \in \mathbb{R}^{n \times d}, \hat{c}_t \in \mathbb{R}^n$, and $s_t \in \mathbb{R}^{\left|\mathcal{N}_t\right|}$. Furthermore, $\mathcal{K}$ is a convex set to be determined. This optimization problem essentially characterizes the extent to which the optimality condition \Cref{eq:MOM} is satisfied, that is, it aims to achieve the largest possible margin for satisfying the inequality constraints.
        However, the MOM method is invalid when 
        the observed solution is not at a vertex of the convex feasible region (for example, in the Noisy Decision setting, the observed decision is not even feasible), and therefore the basis of the observed solution cannot be identified.

\clearpage

\section{Additional Experimental Details}\label{sec: Experiment-details-appendix}
In \Cref{sec:experiments}, we provide synthetic data experimental results and real data experimental results for different methods. In this section, we further explain the details of these experiments setup, implementation, and provide additional experimental results. All experiments are implemented on a cloud computing platform with 128 CPUs of model Intel(R) Xeon(R) Platinum 8369B CPU @ 2.70GHz, 250GB RAM and 500GB storage. 

\subsection{Synthetic Data Experiment Setups and Implementation details} \label{sub:synthetic-experiment}
We present several numerical results, including both linear and nonlinear examples, that demonstrate the great advantages in time consuming, parameter estimation accuracy, optimal solution gap and objective function value gap. Specifically, we study three inverse optimization scenarios for each example.
\begin{itemize}
    \item \textbf{Noisy Decision}: The observed solution is equal to the real optimal solution plus a noise drawing from the standard normal distribution. For example, if true optimization solution for the forward optimization problem(FOP) is $x^* \in \mathbb{R}^p$, then the observed solution is $y=x^* + w$, where every element of $w \in \mathbb{R}^p$ draws from a standard normal distribution.  It is worth noting that the observed solution in this case may not even be feasible, so many classical inverse optimization methods do not handle this situation well.
    \item \textbf{Noisy Objective Function}: The observed solution is the optimal solution corresponding to the key term $\theta+u$ of the objective function plus the standard normal distribution noise. For example, if the true forward optimization problem(FOP) is $ \min \{(\theta + u)^{\top} x \mid x \in [0, 1]^{p} \}$, then the observed solution corresponds to the optimal solution of the optimization problem $ \min \{(\theta + u + w )^{\top} x \mid x \in [0, 1]^{p} \}$, where every element of $w \in \mathbb{R}^p$ drawn from a standard normal distribution.
    \item \textbf{Noiseless}: The observed solution is the optimal solution to the true forward optimization problem without any noise.
\end{itemize}
We repeated this 100 times for each example and uniformly used a random number seed with the same number to ensure fairness. We implement various methods with sample sizes of 50, 100, 300, 500, and 1000, utilizing the Gurobi solver for all experiments. The following metrics are collected:
\begin{itemize}
    \item \textbf{Time}: The total time required to complete 100 serial repetitions of the experiment using the specified method.
    
    \item \textbf{Parameter Error}: The $l_1$ norm of the difference between the estimated parameter $\hat{\theta}$ and the true parameter $\theta^*$, denoted as  $E(\hat{\theta}, \theta^{\star}) := \Vert \hat{\theta}-\theta^* \Vert_1$ in \Cref{eq:parameter-error}.
    
    \item \textbf{Decision Error}: The squared $l_2$ norm of the difference between the optimal solution $x^{\star}(\hat{\theta}; u)$ based on the estimated $\hat{\theta}$ and the optimal solution $x^{\star}(\theta^{\star}; u)$ based on the true parameter $\theta^*$ given context $u$, denoted $D(\hat{\theta}, \theta^{\star}; u) := \Vert x^{\star}(\hat{\theta}; u) - x^{\star}(\theta^{\star}; u) \Vert_2^2$ as shown in \Cref{sub:decision-error}.
    
    \item \textbf{Regret}: The gap between the objective function value of the actual decision and the true optimal objective function value. For example, given context $u$, if the forward optimization problem is $\min \{ h(\theta^*; u)^{\top} x \mid g(u, x) \preceq 0 \} = \max \{ -h(\theta^*; u)^{\top} x \mid g(u, x) \preceq 0 \}$ with the corresponding optimization solution $x^{\star}(\theta^{\star}; u)$ 
    the conditional decision regret is then $\operatorname{Reg}(\hat{\theta}, \theta^{\star}; u) := h(\theta^*; u)^{\top} x^{\star}(\hat{\theta}; u) - h(\theta^*; u)^{\top} x^{\star}(\theta^{\star}; u) $ as shown in \Cref{eq:decision-regret}.
    
\end{itemize}
To fully illustrate the effectiveness and advantages of our method, in each set of experiments, we compared the Fenchel-Young loss method with a number of classical inverse optimization methods, which are SPA \cite{aswani2018inverse}, approximate KKT conditions(KKA) \cite{keshavarz2011imputing}, variational inequality approximate(VIA) \cite{bertsimas2015data}, and maximum optimality margin(MOM) \cite{sun2023maximum}. Among these, the maximum optimality margin (MOM) method models the basic and non-basic variables in the observed solution when the forward optimization problem (FOP) is a linear programming problem. However, in the Noisy Decision setting, the observed optimal solution is not at the vertex of the feasible set. Therefore, this method is not suitable for the Noisy Decision setting, nor is it suitable when the FOP is a nonlinear programming problem.

We design five different synthetic data experiments. The first three experiments involve linear programming forward optimization problems, which include box constraints and half-space constraints. The objective functions are of the form $\theta+u$ and $\theta \circ u$, where $\theta$ is the unknown parameter vector to be estimated, and $u$ is the known context vector. The last two experiments involve forward optimization problems of nonlinear programming: the first is quadratic in $x$ with a box constraint, and the second is linear in $X$ with a feasible set defined by a ball constraint.

We will show that our method is the best, or close to the best, in Noisy Decision setting, Noisy Objective Function setting, or Noiseless setting. Compared with other classical inverse optimization methods, in all kinds of forward optimization problem(FOP), it has significant advantages in parameter recovery, decision making and computational time.

\subsubsection{Example A}
In this experiment, we consider optimization problems with more complex constraints, 
$$
    \min \{ (\theta + u)^{\top} x \mid x \succeq 0, \parallel x \parallel_1 \leq a \} = \min \left\{\sum_{k=1}^{p} (\theta_k + u_k) x_k \mid x \succeq 0, \parallel x \parallel_1 \leq a \right\}.
$$
Every element of context vector $u \in \mathbb{R}^p$ has a uniform distribution with support $[-1,1]^p$ and $p=10$. In addition, we set $a=3$ and true parameter vector $\theta=[0.5,\cdots, 0.5] \in \mathbb{R}^p$. The optimal solution $x^{\text{opt}}$ to this problem is
$$
    x^{\star}(\theta; u)_k = \left\{
    \begin{aligned}
        & a, && \text{if } (\theta_k + u_k) \leq 0 \\
        & 0, && \text{otherwise}.
    \end{aligned}
    \right.
$$
To compute the Fenchel-Young loss, we first derive the convex conjugate of $\Omega$, which is
$$
    \max_{x \in \mathcal{X}(u)} V_{\lambda}(\theta; u, x) = \max_{x \in \mathcal{F} } \left\{ -(\theta + u)^{\top} x - \frac{\lambda}{2} x^{\top}x \right\} ~~~ \rightarrow ~~~ x_{\lambda}^{\star}(\theta; u).
$$
Here $\mathcal{F}:\{x \in \mathbb{R}^p| x \succeq 0, \sum_{i=1}^p x_i \leq a \}$ is the feasible region of this example, and $x_{\lambda}^{\star}(\theta; u)$ is the regularized optimal solution. 
Therefore, based on the observation solution $y$, we have the following expression of Fenchel-Young loss,
$$
\begin{aligned}
    L_{\lambda}(\theta;u,y) &= \max_{x \in \mathcal{X}(u)} V_{\lambda}(\theta; u, x) + \lambda\Omega(y) - \left[ - \sum_{k=1}^{p} (\theta_k + u_k) y_{k} \right] \\ 
    &= \max_{x \in \mathcal{X}(u)} V_{\lambda}(\theta; u, x) + \lambda\Omega(y) + \sum_{k=1}^{p} (\theta_k + u_k) y_{k}
\end{aligned}
$$
Then, we have the gradient of Fenchel-Young Loss:
$$
    \nabla_{\theta_k} L_{\lambda}(\theta;u,y) =  \big[- x_{\lambda}^{\star}(\theta; u)_k + y_k \big]
$$
for $k = 1, \dots, p$. In this experiment, we cannot directly write the explicit expression for $x^*(\theta,u)$, so we need to deploy the solver every time in the experiment, which may result in a long computational time. However, we can choose to use a stochastic gradient descent algorithm to speed up our method and still ensure that the computational time of our method is acceptable.

\paragraph{Result Discussion} The results of experiments in the Noisy Decision setting, the Noisy Objective Function setting, and the Noiseless setting, which are shown in \Cref{tab:example-A-Noisy-Solution-appendix}, \Cref{tab:example-A-Noisy-objective-function-appendix} and \Cref{tab:example-A-Noiseless-appendix} respectively, average over 100 repetitions. Based on the experimental results, we can draw the following conclusions.
\begin{itemize}
    \item In the Noisy Decision setting and the Noisy Objective Function setting, our method based on the Fenchel-Young loss demonstrates superior performance in terms of parameter error, decision error, and regret. Specifically, it achieves the lowest values for these metrics compared to other traditional inverse optimization methods. Other methods, except for KKA, show gradual convergence in decision error and regret as the sample size increases. However, their performance in parameter recovery remains consistently poor. In contrast, the KKA method fails to converge in terms of parameter error, decision error, and regret. It seems that KKA always performs well in the noiseless setting and poorly in the noisy setting in the following experiments.
    \item In the Noiseless setting, the results of our method are very close to those of KKA and VIA, which are considered the best-performing methods. Additionally, it appears that MOM is better suited for decision-making than for parameter recovery. Furthermore, the SPA method performs poorly in terms of parameter error, decision error, and regret. This is because the semi-parametric algorithm in SPA is designed to remove noise. However, in the Noiseless setting, it becomes detrimental to the algorithm's performance.
    
\end{itemize}

\begin{table}[ht]
\caption{The experimental results for various methods under different metrics are presented for the FOP defined as $\min \{ (\theta + u)^{\top} x \mid x \succeq 0, \parallel x \parallel_1 \leq a \} = \min \left\{\sum_{k=1}^{p} (\theta_k + u_k) x_k \mid x \succeq 0, \parallel x \parallel_1 \leq a \right\}$ in the \textbf{Noisy Decision} setting. For the FY loss, we set $\lambda=0.1$ and $\Omega(x) = 1/2\|x\|_2^2$.}
\vskip 0.15in
\centering
\begin{tabular}{c|cccc|cccc|cccc}
\hline
\multirow{2}{*}{Sample size} & \multicolumn{4}{c|}{Parameter Error} & \multicolumn{4}{c|}{Decision Error} & \multicolumn{4}{c}{Regret} \\ \cline{2-13} 
 & FY & SPA & KKA & VIA & FY & SPA & KKA & VIA & FY & SPA & KKA & VIA \\ \hline
50 & \textbf{1.30} & 2.46 & 4.25 & 3.81 & \textbf{3.03} & 4.53 & 7.66 & 3.90 & \textbf{0.05} & 0.20 & 0.77 & 0.08 \\
100 & \textbf{1.17} & 1.35 & 4.58 & 4.12 & \textbf{2.47} & 2.94 & 7.79 & 3.17 & \textbf{0.03} & 0.06 & 0.84 & 0.05 \\
300 & \textbf{0.53} & 2.10 & 4.87 & 4.42 & \textbf{2.05} & 2.66 & 8.22 & 2.19 & \textbf{0.02} & 0.08 & 0.93 & 0.03 \\
500 & \textbf{0.53} & 2.44 & 4.92 & 4.52 & \textbf{1.64} & 2.69 & 8.32 & 1.88 & \textbf{0.01} & 0.10 & 0.95 & 0.03 \\
1000 & \textbf{0.38} & 2.59 & 4.96 & 4.56 & \textbf{1.33} & 2.64 & 8.40 & 1.60 & \textbf{0.01} & 0.11 & 0.96 & 0.03 \\ \hline
\end{tabular}
\label{tab:example-A-Noisy-Solution-appendix}
\end{table}

\begin{table}[ht]
\caption{The experimental results for various methods under different metrics are presented for the FOP defined as $\min \{ (\theta + u)^{\top} x \mid x \succeq 0, \parallel x \parallel_1 \leq a \} = \min \left\{\sum_{k=1}^{p} (\theta_k + u_k) x_k \mid x \succeq 0, \parallel x \parallel_1 \leq a \right\}$ in the \textbf{Noisy Objective Function} setting. For the FY loss, we set $\lambda=0.1$ and $\Omega = 1/2\|x\|_2^2$.}
\vskip 0.15in
\centering
\resizebox{\textwidth}{14mm}{
\begin{tabular}{c|ccccc|ccccc|ccccc}
\hline
\multirow{2}{*}{Sample size} & \multicolumn{5}{c|}{Parameter Error} & \multicolumn{5}{c|}{Decision Error} & \multicolumn{5}{c}{Regret} \\ \cline{2-16} 
 & FY & SPA & KKA & VIA & MOM & FY & SPA & KKA & VIA & MOM & FY & SPA & KKA & VIA & MOM \\ \hline
50 & \textbf{1.15} & 1.84 & 2.17 & 2.09 & 5.91 & \textbf{3.75} & 5.13 & 5.50 & 5.47 & 6.53 & \textbf{0.04} & 0.15 & 0.15 & 0.11 & 0.17 \\
100 & \textbf{1.09} & 0.99 & 1.93 & 1.91 & 5.81 & \textbf{2.81} & 3.63 & 4.46 & 4.00 & 5.37 & \textbf{0.04} & 0.05 & 0.11 & 0.06 & 0.11 \\
300 & \textbf{0.37} & 1.15 & 2.23 & 1.84 & 5.83 & \textbf{2.06} & 2.32 & 4.21 & 2.49 & 3.55 & \textbf{0.03} & 0.03 & 0.14 & 0.03 & 0.06 \\
500 & \textbf{0.36} & 1.43 & 2.89 & 1.69 & 5.81 & \textbf{1.66} & 2.08 & 4.70 & 2.06 & 3.04 & \textbf{0.01} & 0.03 & 0.21 & 0.02 & 0.05 \\
1000 & \textbf{0.22} & 1.48 & 3.50 & 1.81 & 5.82 & \textbf{1.25} & 1.71 & 5.46 & 1.61 & 2.38 & \textbf{0.00} & 0.03 & 0.33 & 0.02 & 0.04 \\ \hline
\end{tabular}}
\label{tab:example-A-Noisy-objective-function-appendix}
\end{table}

\begin{table}[ht]
\caption{The experimental results for various methods under different metrics are presented for the FOP defined as $\min \{ (\theta + u)^{\top} x \mid x \succeq 0, \parallel x \parallel_1 \leq a \} = \min \left\{\sum_{k=1}^{p} (\theta_k + u_k) x_k \mid x \succeq 0, \parallel x \parallel_1 \leq a \right\}$ in the \textbf{Noiseless} setting. For the FY loss, we set $\lambda=0.1$ and $\Omega(x) = 1/2 \|x\|_2^2$.}
\vskip 0.15in
\centering
\resizebox{\textwidth}{14mm}{
\begin{tabular}{c|ccccc|ccccc|ccccc}
\hline
\multirow{2}{*}{Samplesize} & \multicolumn{5}{c|}{ParameterError} & \multicolumn{5}{c|}{DecisionError} & \multicolumn{5}{c}{Regret} \\ \cline{2-16} 
 & FY & SPA & KKA & VIA & MOM & FY & SPA & KKA & VIA & MOM & FY & SPA & KKA & VIA & MOM \\ \hline
50 & \textbf{0.64} & 2.12 & 0.96 & 0.93 & 5.01 & \textbf{1.05} & 3.02 & 0.00 & 0.00 & 3.12 & \textbf{0.02} & 0.10 & 0.00 & 0.00 & 0.06 \\
100 & \textbf{0.70} & 2.33 & 0.44 & 0.42 & 5.03 & \textbf{1.09} & 2.89 & 0.00 & 0.00 & 2.76 & \textbf{0.01} & 0.10 & 0.00 & 0.00 & 0.05 \\
300 & \textbf{0.34} & 2.67 & 0.16 & 0.15 & 5.07 & \textbf{0.76} & 2.87 & 0.00 & 0.00 & 2.07 & \textbf{0.00} & 0.12 & 0.00 & 0.00 & 0.03 \\
500 & \textbf{0.21} & 2.69 & 0.08 & 0.07 & 5.08 & \textbf{0.63} & 2.72 & 0.00 & 0.00 & 1.73 & \textbf{0.00} & 0.12 & 0.00 & 0.00 & 0.03 \\
1000 & \textbf{0.18} & 2.71 & 0.04 & 0.04 & 5.08 & \textbf{0.53} & 2.56 & 0.00 & 0.00 & 1.40 & \textbf{0.00} & 0.12 & 0.00 & 0.00 & 0.03 \\ \hline
\end{tabular}}
\label{tab:example-A-Noiseless-appendix}
\end{table}

\clearpage

\subsubsection{Example B}  \label{sec:example-B-details}
Next, we consider the following forward optimization problem(FOP),
$$
    \min \{ (\theta \circ u)^{\top} x \mid x \in [-1, 1]^{p} \} = \min \left\{\sum_{k=1}^{p} (\theta_k \circ u_k) x_k \mid x \in [-1, 1]^{p} \right\}.
$$

Every element of context vector $u \in \mathbb{R}^p$ has a uniform distribution with support $[-1,1]^p$ and $p=10$ in our experiments. $\circ$ means Hadamard product, that is, $(\theta \circ u)^{\top} x = \sum_{i=1}^p \theta_i u_i x_i$. The optimal solution $x_{\text{opt}}$ to this problem is
$$
    x^{\star}(\theta; u)_k = \left\{
    \begin{aligned}
        & -1, && \text{if } \theta_ku_k > 0 \\
        & +1, && \text{if } \theta_ku_k < 0.
    \end{aligned}
    \right.
$$
Note that this problem is not identifiable for the parameter vector $\theta \in \mathbb{R}^p$ that needs to be estimated. Because, as long as we estimate the sign of each component of the parameter vector $\theta$ correctly, we can obtain the solution that is exactly consistent with the optimal parameter, which also means the parameter vector $\theta$ is not identifiable. Therefore, we set the true parameter vector $\theta=[0.5,-0.5,\cdots,-0.5,0.5,1] \in \mathbb{R}^p$ to test the sign recognition ability of various inverse optimization methods. 

The last element is set to 1 instead of 0.5 because the parameter vector $\theta$ solved by the SPA, KKA, and VIA methods will be very close to 0, which may cause the degeneracy. The reasons for this will be explained in the next few paragraphs. Therefore, we need to ensure that the sum of the elements of the real parameter vector $\theta$ is not 0, so that we can add constraints such as the sum of the elements of $\theta$ is a non-zero constant to prevent the degeneracy of these methods.

To compute the Fenchel-Young loss, we first derive the convex conjugate of $\Omega$, which is
$$
    \max_{x \in \mathcal{X}(u)} V_{\lambda}(\theta; u, x) = \max_{x \in [-1,1]^{p} } \left\{ -(\theta \circ u)^{\top} x - \frac{\lambda}{2} x^{\top}x \right\} ~~~ \rightarrow ~~~ x_{\lambda}^{\star}(\theta; u)_k = \max \left\{-1,  \min \left\{1, - \frac{\theta_k u_k}{\lambda} \right\} \right\}.
$$
Therefore, based on the observation solution $y$, we have the following expression of Fenchel-Young loss,
$$
\begin{aligned}
    L_{\lambda}(\theta;u,y) &= \max_{x \in \mathcal{X}(u)} V_{\lambda}(\theta; u, x) + \lambda\Omega(y) - \left[ - \sum_{k=1}^{p} (\theta_k \circ u_k) y_{k} \right] \\ 
    &= \max_{x \in \mathcal{X}(u)} V_{\lambda}(\theta; u, x) + \lambda\Omega(y) + \sum_{k=1}^{p} (\theta_k \circ u_k) y_{k}
\end{aligned}
$$
Then, we have the gradient of Fenchel-Young Loss: for $k = 1,\cdots,p$,
$$
    \nabla_{\theta_k} L_{\lambda}(\theta;u,y) =  \big[- u_kx_{\lambda}^{\star}(\theta; u)_k + u_ky_k \big] = u_k\left[- \max \left\{-1,  \min \left\{1, -\frac{\theta_k u_k}{\lambda} \right\} \right\} + y_k \right].
$$

Before showing the experimental results of this example, let us first explain why the KKA, VIA, and SPA methods degenerate. The following is the Lagrangian dual of the forward optimization problem(FOP).
$$
\begin{aligned}
    L(x, \lambda_1, \lambda_2, u, \theta) &= (\theta \circ u)^{\top} x + \lambda_1^{\top}(-1 - x) + \lambda_2^{\top}(x - 1) \\
    &= (\theta \circ u - \lambda_1 + \lambda_2)^{\top} x - \sum_{k=1}^{p} (\lambda_{1,k} + \lambda_{2,k}) \\
    &= \sum_{k=1}^{p} (\theta \circ u - \lambda_1 + \lambda_2)_k x_k - \sum_{k=1}^{p} (\lambda_{1,k} + \lambda_{2,k}).
\end{aligned}
$$
where $\lambda_1$ and $\lambda_2$ are Lagrangian multipliers. By minimizing the $x$, we obtain the Lagrangian dual function as
$$
  h(\lambda_1, \lambda_2, u, \theta) = \inf_{x} L(x, \lambda_1, \lambda_2, u, \theta) =
  \left\{
  \begin{aligned}
    & - \sum_{p=1}^{k} (\lambda_{1,k} + \lambda_{2,k}), && (\theta \circ u - \lambda_1 + \lambda_2)_k = 0, \ \forall k \in [p] \\
    & - \infty, && \text{otherwise}.
  \end{aligned}
  \right.
$$
\paragraph{SPA Method Formulation} Therefore, for the SPA method, we are dealing with the following constrained risk minimization problem: 
$$
\begin{aligned}
    \hat{\theta}_n = \arg\min_{\theta, \lambda_{i,1}, \lambda_{i,2}, \epsilon_i} &\frac{1}{n} \sum_{i=1}^n \epsilon_i \\
    \text { s.t. } &  \sum_{k=1}^{p} (\theta  u_i)_k \bar{y}_{i,k} - \Big[-\sum_{k=1}^{p} (\lambda_{i,1} + \lambda_{i,2})_k \Big] \leq \epsilon_i , && \forall i \in [n] \\
    & (\theta \circ u_i - \lambda_{i, 1} + \lambda_{i, 2})_k = 0, && \forall k \in [p], ~ \forall i \in [n] \\
    & \lambda_{i,1} \succeq 0, && \forall i \in [n] \\
    & \lambda_{i,2} \succeq 0, && \forall i \in [n] 
    \label{eq:multiply-SPA}
\end{aligned}
$$
where $\bar{y}_{i,k}$ is the value of the $k$-th element of the $i$-th data observation solution vector after noise removal by the semi-parametric algorithm.

\paragraph{VIA Method Formulation} For the VIA method, given a sequence of observation $y_j$ for $j = 1, \dots, N$, we have the following optimization problem:
$$
\begin{aligned}
    \min_{\theta, \lambda_{1,i}, \lambda_{2,i}, \epsilon_i } ~~ & \|\epsilon \|_2^2 \\
    \text{s.t.} ~~ & (\theta \circ u_i)^{\top} y_i - \left[- \sum_{k=1}^{p} (\lambda_{1,i,k} + \lambda_{2,i,k}) \right] \leq \epsilon_i, ~~ && \forall i \in [n] \\
    & (\theta \circ u_i - \lambda_{1, i} + \lambda_{2, i})_k = 0, ~~ && \forall k \in [p], ~ \forall i \in [n] \\
    & \lambda_{1,j} \succeq 0, \lambda_{2, j} \succeq 0 ~~ && \forall i \in [n] 
    \label{eq:multiply-VIA}
\end{aligned}
$$

\paragraph{KKA Method Formulation} For the KKA method, we need to solve the following optimization problem.
$$
\begin{aligned}
    \min_{\theta, \lambda_{1,i}, \lambda_{2,i} } ~~ & \sum_{i=1}^{n} \big( \|r_i^{\rm stat} \|_2^2 + \|r_{1,i}^{\rm comp}\|_2^2 + \|r_{2,i}^{\rm comp}\|_2^2 \big) \\
    \text{s.t.} ~~ & r_{i, k}^{\rm stat}(\theta, \lambda_i, \nu_i) = \theta_ku_{i,k} + \lambda_{1, i} - \lambda_{2, i}, && \forall k \in [p], ~ \forall i \in [n]\\
    & r_{1,i,k}^{\rm comp}(\lambda_i ) = \lambda_{1,i,k} (x_{i,k} - 1), && \forall k \in [p], ~ \forall i \in [n] \\
    & r_{2,i,k}^{\rm comp}(\lambda_i ) = \lambda_{2,i,k} (-1 - x_{i,k}), && \forall k \in [p], ~ \forall i \in [n] \\
    & \lambda_i \succeq 0, && \forall i \in [n] \\
    & \theta \in [0, 1]^{p}.
\end{aligned}
$$

Based on the above optimization problems, we can find that when $\theta$ is equal to a vector with all zeros, the objective function of these three optimization problems can be taken to 0, that is, the minimum value. As a result, this can mislead the solver into solving meaningless solutions such as all zeros, leading to degeneracy of the problem.

\paragraph{Result Discussion} The results of experiments in the Noisy Decision setting, the Noisy Objective Function setting, and the Noiseless setting, which are shown in \Cref{tab:example-B-Noisy-Solution-appendix}, \Cref{tab:example-B-Noisy-objective-function-appendix} and \Cref{tab:example-B-Noiseless-appendix} respectively, average over 100 repetitions. Since the parameter vector $\theta$ is not identifiable in this problem, we only need to focus on the decision error and the regret. Thus, we can clearly come to the following conclusion.

\begin{itemize}
    \item In the Noisy Decision setting and Noisy Objective Function setting, the Fenchel-Young loss and VIA methods exhibit the best performance, achieving zero decision error and regret even with small sample sizes. In contrast, the MOM method requires a larger sample size to reduce these metrics to zero. The KKA and SPA methods perform the worst, with almost no convergence observed.
       
    Furthermore, we observe that the parameter $\theta$ estimated by the SPA, KKA, and VIA methods consistently approaches zero, leading to the problem degeneracy. The reason for the degeneracy, as previously explained, is that $\theta=\textbf{0}$ happens to minimize the optimization objectives of SPA, KKA, and VIA. However, such a $\theta$ is meaningless for decision-making. Among these methods, VIA is still able to correctly recover the sign of each element of the vector $\theta$, which helps to maintain good performance in terms of decision error and regret.
    
    To address this issue, we attempted to introduce regularization constraints to these methods. For example, we imposed constraints such as requiring the sum of the elements of the estimated parameter vector $\theta$ to be 0.5, or fixing the first element of the parameter vector $\theta$ at 0.5. However, these modifications did not significantly improve the experimental results of these methods. In contrast, our proposed method based on the Fenchel-Young loss does not face the challenge of degeneracy and maintains robust performance without the need for such regularization constraints. Notably, MOM also successfully avoids degeneration, which is in agreement with the conclusions reached by \citet{sun2023maximum}.
    
    \item  In the Noiseless setting, all methods except SPA achieve very good results, with decision error and regret reaching zero even when the sample size is as small as 50. However, SPA incorporates a semi-parametric algorithm designed to remove noise, which becomes counterproductive in a noiseless environment. If the semi-parametric algorithm is removed, the method essentially becomes equivalent to VIA. Additionally, compared to the noisy settings, KKA is able to obtain a correctly signed parameter vector instead of a degenerate all-zero vector in this noiseless setting.
    
\end{itemize}



\begin{table}[ht]
\caption{The experimental results for various methods under different metrics are presented for the FOP defined as $ \min \{ (\theta \circ u)^{\top} x \mid x \in [-1, 1]^{p} \}= \min \left\{\sum_{k=1}^{p} (\theta_k * u_k) x_k \mid x \in [-1, 1]^{p} \right\}$ in the \textbf{Noisy Decision} setting. For the FY loss, we set $\lambda=0.1$ and $\Omega(x) = 1/2\|x\|_2^2$.}
\vskip 0.15in
\centering
\begin{tabular}{c|cccc|cccc|cccc}
\hline
\multirow{2}{*}{Sample size} & \multicolumn{4}{c|}{Parameter Error} & \multicolumn{4}{c|}{Decision Error} & \multicolumn{4}{c}{Regret} \\ \cline{2-13} 
 & FY & SPA & KKA & VIA & FY & SPA & KKA & VIA & FY & SPA & KKA & VIA \\ \hline
50 & \textbf{2.76} & 5.50 & 5.50 & 5.50 & \textbf{0.00} & 19.99 & 19.99 & 0.00 & \textbf{0.00} & 2.75 & 2.75 & 0.00 \\
100 & \textbf{2.42} & 5.50 & 5.50 & 5.50 & \textbf{0.00} & 19.99 & 19.99 & 0.00 & \textbf{0.00} & 2.76 & 2.76 & 0.00 \\
300 & \textbf{2.25} & 5.50 & 5.50 & 5.50 & \textbf{0.00} & 20.00 & 20.00 & 0.00 & \textbf{0.00} & 2.75 & 2.75 & 0.00 \\
500 & \textbf{2.27} & 5.50 & 5.50 & 5.50 & \textbf{0.00} & 20.01 & 20.01 & 0.00 & \textbf{0.00} & 2.76 & 2.76 & 0.00 \\
1000 & \textbf{2.19} & 5.50 & 5.50 & 5.50 & \textbf{0.00} & 19.99 & 19.99 & 0.00 & \textbf{0.00} & 2.75 & 2.75 & 0.00 \\ \hline
\end{tabular}
\label{tab:example-B-Noisy-Solution-appendix}
\end{table}

\begin{table}[ht]
\caption{The experimental results for various methods under different metrics are presented for the FOP defined as $ \min \{ (\theta \circ u)^{\top} x \mid x \in [-1, 1]^{p} \}= \min \left\{\sum_{k=1}^{p} (\theta_k * u_k) x_k \mid x \in [-1, 1]^{p} \right\}$ in the \textbf{Noisy Objective Function} setting. For the FY loss, we set $\lambda = 0.1$ and $\Omega = 1/2\|x\|_2^2$.}
\vskip 0.15in
\centering
\resizebox{\textwidth}{14mm}{
\begin{tabular}{c|ccccc|ccccc|ccccc}
\hline
\multirow{2}{*}{Sample size} & \multicolumn{5}{c|}{Parameter Error} & \multicolumn{5}{c|}{Decision Error} & \multicolumn{5}{c}{Regret} \\ \cline{2-16} 
 & FY & SPA & KKA & VIA & MOM & FY & SPA & KKA & VIA & MOM & FY & SPA & KKA & VIA & MOM \\ \hline
50 & \textbf{2.78} & 5.50 & 5.50 & 5.50 & 7.94 & \textbf{0.00} & 19.99 & 19.99 & 2.04 & 1.96 & \textbf{0.00} & 2.75 & 2.75 & 0.25 & 0.24 \\
100 & \textbf{2.80} & 5.50 & 5.50 & 5.50 & 7.27 & \textbf{0.00} & 19.99 & 19.99 & 0.20 & 0.16 & \textbf{0.00} & 2.76 & 2.76 & 0.03 & 0.02 \\
300 & \textbf{2.80} & 5.50 & 5.50 & 5.50 & 7.22 & \textbf{0.00} & 20.00 & 20.00 & 0.00 & 0.00 & \textbf{0.00} & 2.75 & 2.75 & 0.00 & 0.00 \\
500 & \textbf{2.80} & 5.50 & 5.50 & 5.50 & 7.22 & \textbf{0.00} & 20.01 & 20.01 & 0.00 & 0.00 & \textbf{0.00} & 2.76 & 2.76 & 0.00 & 0.00 \\
1000 & \textbf{2.80} & 5.50 & 5.50 & 5.50 & 7.17 & \textbf{0.00} & 19.99 & 19.99 & 0.00 & 0.00 & \textbf{0.00} & 2.75 & 2.75 & 0.00 & 0.00 \\ \hline
\end{tabular}}
\label{tab:example-B-Noisy-objective-function-appendix}
\end{table}

\begin{table}[ht]
\caption{The experimental results for various methods under different metrics are presented for the FOP defined as $ \min \{ (\theta \circ u)^{\top} x \mid x \in [-1, 1]^{p} \}= \min \left\{\sum_{k=1}^{p} (\theta_k * u_k) x_k \mid x \in [-1, 1]^{p} \right\}$ in the \textbf{Noiseless} setting. For the FY loss, we set $\lambda=0.1$ and $\Omega(x) = 1/2\|x\|_2^2$.}
\vskip 0.15in
\centering
\resizebox{\textwidth}{14mm}{
\begin{tabular}{c|ccccc|ccccc|ccccc}
\hline
\multirow{2}{*}{Sample size} & \multicolumn{5}{c|}{Parameter Error} & \multicolumn{5}{c|}{Decision Error} & \multicolumn{5}{c}{Regret} \\ \cline{2-16} 
 & FY & SPA & KKA & VIA & MOM & FY & SPA & KKA & VIA & MOM & FY & SPA & KKA & VIA & MOM \\ \hline
50 & \textbf{2.22} & 5.50 & 4.50 & 4.33 & 31.27 & \textbf{0.00} & 20.00 & 0.00 & 0.00 & 0.00 & \textbf{0.00} & 20.00 & 0.00 & 0.00 & 0.00 \\
100 & \textbf{2.22} & 5.50 & 4.50 & 4.39 & 31.24 & \textbf{0.00} & 20.00 & 0.00 & 0.00 & 0.00 & \textbf{0.00} & 20.00 & 0.00 & 0.00 & 0.00 \\
300 & \textbf{2.22} & 5.50 & 4.50 & 4.45 & 31.26 & \textbf{0.00} & 20.00 & 0.00 & 0.00 & 0.00 & \textbf{0.00} & 20.00 & 0.00 & 0.00 & 0.00 \\
500 & \textbf{2.22} & 5.50 & 4.50 & 4.48 & 31.29 & \textbf{0.00} & 20.01 & 0.00 & 0.00 & 0.00 & \textbf{0.00} & 20.01 & 0.00 & 0.00 & 0.00 \\
1000 & \textbf{2.22} & 5.50 & 4.50 & 4.49 & 31.32 & \textbf{0.00} & 20.00 & 0.00 & 0.00 & 0.00 & \textbf{0.00} & 20.00 & 0.00 & 0.00 & 0.00 \\ \hline
\end{tabular}}
\label{tab:example-B-Noiseless-appendix}
\end{table}
\clearpage
\subsubsection{Example C} 
Next, we consider the following forward optimization problem (FOP), 
\begin{align*}
    \min \{ (\theta + u)^{\top} x \mid x \in [-1, 1]^{p} \} = \min \left\{\sum_{k=1}^{p} (\theta_k + u_k) x_k \mid x \in [-1, 1]^{p} \right\}.
\end{align*}
Every element of context vector $u \in \mathbb{R}^p$ has a uniform distribution with support $[-1,1]^p$ and $p=10$ in our experiments. Apparently, the optimal solution $x^{\text{opt}}$ to this problem is
$$
    x^{\star}(\theta; u)_k = \left\{
    \begin{aligned}
        & -1, && \text{if } \theta_k + u_k > 0 \\
        & +1, && \text{if } \theta_k + u_k \leq 0.
    \end{aligned}
    \right.
$$
 We set the true parameter vector $\theta=[0.5, \cdots, 0.5]^p \in \mathbb{R}^p$, which needs to be estimated. In addition, our method uniformly uses the $l_2$ regular term, which is $\frac{\lambda}{2}\|x\|_2^2$ , in all experiments.

To compute the Fenchel-Young loss, we first derive the convex conjugate of $\Omega$, which is
$$
    \max_{x \in \mathcal{X}(u)} V_{\lambda}(\theta; u, x) = \max_{x \in [-1,1]^{p} } \left\{ -(\theta + u)^{\top} x - \frac{\lambda}{2} x^{\top}x \right\} ~~~ \rightarrow ~~~ x_{\lambda}^{\star}(\theta; u)_k = \max \left\{-1,  \min \left\{1, - \frac{\theta_k + u_k}{\lambda} \right\} \right\}.
$$

Therefore, based on the observation solution $y$, we have the following expression of Fenchel-Young loss,
$$
    L_{\lambda}(\theta;u,y) = \max_{x \in \mathcal{X}(u)} V_{\lambda}(\theta; u, x) + \lambda\Omega(y) - \left[ - \sum_{k=1}^{p} (\theta_k + u_k) y_{k} \right] = \max_{x \in \mathcal{X}(u)} V_{\lambda}(\theta; u, x) + \lambda\Omega(y) + \sum_{k=1}^{p} (\theta_k + u_k) y_{k}
$$
Then, we have the gradient of Fenchel-Young Loss:
$$
    \nabla_{\theta_k} L_{\lambda}(\theta;u,y) =  \big[- x_{\lambda}^{\star}(\theta; u)_k + y_k \big] = \left[- \max \left\{-1,  \min \left\{1, \frac{\theta_k + u_k}{\lambda} \right\} \right\} + y_k \right]
$$
for $k = 1, \dots, p$. Finally, we can utilize the gradient descent algorithm to optimize the parameter vector $\theta$.

\paragraph{Result Discussion}
The results of experiments in the Noisy Decision setting, the Noisy Objective Function setting, and the Noiseless setting, which are shown in \Cref{tab:example-C-Noisy-Solution-appendix}, \Cref{tab:example-C-Noisy-objective-function-appendix} and \Cref{tab:example-C-Noiseless-appendix} respectively, average over 100 repetitions. Based on these results, we can draw the following conclusions.
\begin{itemize}
    \item In the Noisy Decision setting, our method shows significant advantages over other inverse optimization methods in both parameter estimation and decision-making. The parameter error, decision error, and regret are the lowest among all methods and decrease rapidly as the sample size increases. Moreover, the SPA method outperforms the VIA method in terms of parameter error, decision error, and regret. This superior performance may be attributed to the semi-parametric algorithm in SPA, which effectively mitigates the impact of noise.
    
    \item In the Noisy Objective Function setting, our method achieves the lowest parameter error, decision error, and regret. Specifically, because standard normal distribution noise $w$ is added to $\theta+u$, the problem becomes more complex. If $\theta+u+w \leq 0$, the observed solution is 1; otherwise, it is -1. This scenario is more challenging than the Noisy Decision setting, where the direct addition of normal distribution noise to the observed solution has a relatively small probability of changing 1 to -1 or -1 to 1. In contrast, in the Noisy Objective Function setting, since $u$ is uniformly distributed from $[-1,1]^p$ and $\theta=[0.5,\cdots,0.5]^p$, the addition of standard normal distribution noise $w$ is more likely to alter the original sign of $\theta+u$. Consequently, the noise has a more pronounced impact in the Noisy Objective Function setting.

    Due to these complexities, all methods perform significantly worse in this setting compared to the Noisy Decision setting, and even the KKA method fails to converge. However, our method remains the best performing approach in this challenging scenario.
    
    \item In the Noiseless setting, almost all methods perform well except MOM, particularly in terms of regret, which is expected given its status as a classic setting. However, even with a small sample size, such as 50, our method achieves superior results compared to other classical inverse optimization methods. The inclusion of a semi-parametric algorithm in SPA, designed to remove noise, is redundant in the Noiseless setting and can detract from its performance. In this specific example, removing the semi-parametric algorithm from SPA essentially renders it equivalent to VIA.
\end{itemize}


\begin{table}[ht] 
\caption{The experimental results for various methods under different metrics are presented for the FOP defined as $\min \{ (\theta + u)^{\top} x \mid x \in [-1, 1]^{p} \}=\min \left\{\sum_{k=1}^{p} (\theta_k + u_k) x_k \mid x \in [-1, 1]^{p} \right\}$ in the \textbf{Noisy Decision} setting. For the FY loss, we set $\lambda=0.1$ and $\Omega(x) = 1/2\|x\|_2^2$.}
\vskip 0.15in
\centering
\begin{tabular}{c|cccc|cccc|cccc}
\hline
\multirow{2}{*}{Sample size} & \multicolumn{4}{c|}{ParameterError} & \multicolumn{4}{c|}{DecisionError} & \multicolumn{4}{c}{Regret} \\ \cline{2-13} 
 & FY & SPA & KKA & VIA & FY & SPA & KKA & VIA & FY & SPA & KKA & VIA \\ \hline
50 & \textbf{1.26} & 1.24 & 4.19 & 2.08 & \textbf{2.18} & 2.23 & 9.16 & 3.96 & \textbf{0.11} & 0.12 & 1.07 & 0.29 \\
100 & \textbf{0.82} & 0.86 & 4.61 & 1.88 & \textbf{1.48} & 1.65 & 9.62 & 3.66 & \textbf{0.05} & 0.06 & 1.16 & 0.23 \\
300 & \textbf{0.47} & 0.73 & 4.87 & 1.79 & \textbf{0.87} & 1.46 & 9.86 & 3.58 & \textbf{0.02} & 0.04 & 1.22 & 0.18 \\
500 & \textbf{0.35} & 0.70 & 4.92 & 1.78 & \textbf{0.67} & 1.43 & 9.92 & 3.57 & \textbf{0.01} & 0.03 & 1.23 & 0.17 \\
1000 & \textbf{0.26} & 0.69 & 4.96 & 1.77 & \textbf{0.50} & 1.39 & 9.95 & 3.54 & \textbf{0.01} & 0.03 & 1.24 & 0.16 \\ \hline
\end{tabular}
\label{tab:example-C-Noisy-Solution-appendix}
\end{table}

\begin{table}[h]
\caption{The experimental results for various methods under different metrics are presented for the FOP defined as $\min \{ (\theta + u)^{\top} x \mid x \in [-1, 1]^{p} \}=\min \left\{\sum_{k=1}^{p} (\theta_k + u_k) x_k \mid x \in [-1, 1]^{p} \right\}$ in the \textbf{Noisy Objective Function} setting. For the FY loss, we set $\lambda=0.1$ and $\Omega(x) = 1/2 \|x\|_2^2$.}
\vskip 0.15in
\centering
\resizebox{\textwidth}{14mm}{
\begin{tabular}{c|ccccc|ccccc|ccccc}
\hline
\multirow{2}{*}{Sample size} & \multicolumn{5}{c|}{Parameter Error} & \multicolumn{5}{c|}{Decision Error} & \multicolumn{5}{c}{Regret} \\ \cline{2-16} 
 & FY & SPA & KKA & VIA & MOM & FY & SPA & KKA & VIA & MOM & FY & SPA & KKA & VIA & MOM \\ \hline
50 & \textbf{1.83} & 2.10 & 3.85 & 2.09 & 2.57 & \textbf{3.70} & 4.57 & 7.66 & 4.21 & 5.18 & \textbf{0.25} & 0.35 & 0.78 & 0.30 & 0.38 \\
100 & \textbf{1.67} & 1.94 & 4.38 & 1.99 & 2.47 & \textbf{3.34} & 4.08 & 8.79 & 3.96 & 4.97 & \textbf{0.19} & 0.26 & 0.98 & 0.25 & 0.33 \\
300 & \textbf{1.68} & 2.06 & 4.80 & 2.00 & 2.46 & \textbf{3.39} & 4.21 & 9.57 & 4.02 & 4.96 & \textbf{0.16} & 0.24 & 1.15 & 0.22 & 0.31 \\
500 & \textbf{1.67} & 2.06 & 4.88 & 1.98 & 2.46 & \textbf{3.36} & 4.18 & 9.75 & 3.98 & 4.94 & \textbf{0.15} & 0.23 & 1.19 & 0.21 & 0.31 \\
1000 & \textbf{1.69} & 2.08 & 4.93 & 2.01 & 2.46 & \textbf{3.39} & 4.20 & 9.86 & 4.04 & 4.93 & \textbf{0.15} & 0.23 & 1.22 & 0.21 & 0.31 \\ \hline
\end{tabular}}
\label{tab:example-C-Noisy-objective-function-appendix}
\end{table}

\begin{table}[h]
\caption{The experimental results for various methods under different metrics are presented for the FOP defined as $\min \{ (\theta + u)^{\top} x \mid x \in [-1, 1]^{p} \}=\min \left\{\sum_{k=1}^{p} (\theta_k + u_k) x_k \mid x \in [-1, 1]^{p} \right\}$ in the \textbf{Noiseless} setting. For the FY loss, we set $\lambda=0.01$ and $\Omega(x) = 1/2\|x\|_2^2$.}
\vskip 0.15in
\centering
\resizebox{\textwidth}{14mm}{
\begin{tabular}{c|ccccc|ccccc|ccccc}
\hline
\multirow{2}{*}{Sample size} & \multicolumn{5}{c|}{Parameter Error} & \multicolumn{5}{c|}{Decision Error} & \multicolumn{5}{c}{Regret} \\ \cline{2-16} 
 & FY & SPA & KKA & VIA & MOM & FY & SPA & KKA & VIA & MOM & FY & SPA & KKA & VIA & MOM \\ \hline
50 & \textbf{0.31} & 0.55 & 0.39 & 0.34 & 1.78 & \textbf{0.00} & 1.04 & 0.80 & 0.46 & 3.34 & \textbf{0.00} & 0.03 & 0.02 & 0.07 & 0.16 \\
100 & \textbf{0.15} & 0.62 & 0.20 & 0.19 & 1.74 & \textbf{0.01} & 1.26 & 0.40 & 0.14 & 3.36 & \textbf{0.00} & 0.03 & 0.00 & 0.04 & 0.15 \\
300 & \textbf{0.06} & 0.71 & 0.06 & 0.07 & 1.69 & \textbf{0.03} & 1.45 & 0.13 & 0.04 & 3.37 & \textbf{0.00} & 0.03 & 0.00 & 0.00 & 0.14 \\
500 & \textbf{0.04} & 0.70 & 0.04 & 0.04 & 1.70 & \textbf{0.03} & 1.42 & 0.08 & 0.02 & 3.38 & \textbf{0.00} & 0.03 & 0.00 & 0.00 & 0.14 \\
1000 & \textbf{0.03} & 0.70 & 0.02 & 0.02 & 1.69 & \textbf{0.03} & 1.40 & 0.04 & 0.00 & 3.38 & \textbf{0.00} & 0.03 & 0.00 & 0.00 & 0.14 \\ \hline
\end{tabular}}
\label{tab:example-C-Noiseless-appendix}
\end{table}

\clearpage

\subsubsection{Example D}
Then, we consider a nonlinear programming setting. The forward optimization problem(FOP) is
$$
    \min \{ x^{\top}x - (\theta + u)^{\top} x \mid x \in [0, 1]^{p} \} = \max \left\{-\sum_{k=1}^{p} x_k^2 + \sum_{k=1}^{p} (\theta + u)_k x_k \mid x_k \in [0, 1], \ \forall \ k \in [p] \right\}.
$$
where context vector $u$ follows a uniform distribution with support $[0,2]^{p}$ and $p=10$. We set the true parameter vector $\theta=[0.5,\cdots, 0.5] \in \mathbb{R}^p$. The optimal solution to this problem can be expressed as
$$
x^{\star}(\theta;u)_{k} = \max\left(0, \min\left(1, \frac{ (\theta + u)_k }{2} \right) \right).
$$
To compute the Fenchel-Young loss, we first derive the convex conjugate of $\Omega$, which is
$$
\begin{aligned}
    \max_{x \in \mathcal{X}(u)} V_{\lambda}(\theta; u, x) &= \max_{x \in [0,1]^{p} } \sum_{k=1}^{p} \left\{ -x_k^2 + (\theta_k + u_k)x_k - \frac{\lambda}{2} x_k^2 \right\} ~~~ \rightarrow ~~~ x_{\lambda}^{\star}(\theta;u)_k = \max \left\{0,  \min \left\{1, \frac{\theta_k + u_k}{\lambda + 2} \right\} \right\} 
\end{aligned}
$$
Therefore, based on the observation solution $y$, we have the following expression of Fenchel-Young loss,
$$
\begin{aligned}
    L_{\lambda}(\theta;u,y) &= \max_{x \in \mathcal{X}(u)} V_{\lambda}(\theta; u, x) + \lambda\Omega(y) - \left[ - \sum_{k=1}^{p} (\theta_k + u_k) y_{k} \right] \\ 
    &= \max_{x \in \mathcal{X}(u)} V_{\lambda}(\theta; u, x) + \lambda\Omega(y) + \sum_{k=1}^{p} (\theta_k + u_k) y_{k}
\end{aligned}
$$
Then, we have the gradient of Fenchel-Young Loss:
$$
    \nabla_{\theta_k} L_{\lambda}(\theta;u,y) =  \left[ x_{\lambda}^{\star}(\theta; u)_k - y_k \right] = \left[ \max \left\{0,  \min \left\{1, \frac{\theta_k + u_k}{\lambda + 2} \right\} \right\}  - y_k \right] .
$$
for $k = 1, \dots, p$. 

\paragraph{Result Discussion} The results of experiments in the Noisy Decision setting, the Noisy Objective Function setting, and the Noiseless setting, which are shown in \Cref{tab:example-D-Noisy-Solution-appendix}, \Cref{tab:example-D-Noisy-objective-function-appendix} and \Cref{tab:example-D-Noiseless-appendix} respectively, average over 100 repetitions. Based on the above results, we can draw some interesting conclusions.
\begin{itemize}
    \item In the Noisy Decision and Noisy Objective Function settings, the KKA method even fails to solve when using the Gurobi solver. However, in the Noiseless setting, KKA performs very well. In contrast, our Fenchel-Young loss method consistently outperforms other inverse optimization methods in all settings. This highlights the significant advantages of our method in solving inverse optimization problems, particularly when the forward optimization problem is nonlinear.
    \item In the Noiseless setting, both KKA and our method achieve zero decision error and regret. KKA exhibits the best overall performance, while our method closely follows with a low parameter error. In contrast, VIA performs the worst and does not converge. This may be attributed to the fact that the objective function in this problem is nonlinear with respect to the decision variable $x$. Since VIA employs a first-order conditional approximation, it can only guarantee the equivalence to ASO when the objective function is linear in $x$.
\end{itemize}

\begin{table}[]
\caption{The experimental results for various methods under different metrics are presented for the FOP defined as $\min \{ x^{\top}x - (\theta + u)^{\top} x \mid x \in [0, 1]^{p} \} = \min \left\{\sum_{k=1}^{p} x_k^2 - \sum_{k=1}^{p} (\theta + u)_k x_k \mid x_k \in [0, 1], \ \forall \ k \in [p] \right\}$ in the \textbf{Noisy Decision} setting. For the FY loss, we set $\lambda=0.1$ and $\Omega(x) = 1/2\|x\|_2^2$.}
\vskip 0.15in
\centering
\begin{tabular}{c|cccc|cccc|cccc}
\hline
\multirow{2}{*}{Sample size} & \multicolumn{4}{c|}{Parameter Error} & \multicolumn{4}{c|}{Decision Error} & \multicolumn{4}{c}{Regret} \\ \cline{2-13} 
 & FY & SPA & KKA & VIA & FY & SPA & KKA & VIA & FY & SPA & KKA & VIA \\ \hline
50 & \textbf{2.33} & 2.90 & - & 7.94 & \textbf{0.14} & 0.21 & - & 0.74 & \textbf{0.15} & 0.22 & - & 0.74 \\
100 & \textbf{1.57} & 2.05 & - & 8.01 & \textbf{0.07} & 0.12 & - & 0.75 & \textbf{0.07} & 0.12 & - & 0.75 \\
300 & \textbf{0.92} & 1.52 & - & 7.79 & \textbf{0.02} & 0.07 & - & 0.74 & \textbf{0.02} & 0.07 & - & 0.74 \\
500 & \textbf{0.70} & 1.44 & - & 7.64 & \textbf{0.01} & 0.06 & - & 0.72 & \textbf{0.01} & 0.06 & - & 0.72 \\
1000 & \textbf{0.52} & 1.36 & - & 7.74 & \textbf{0.01} & 0.05 & - & 0.74 & \textbf{0.01} & 0.05 & - & 0.74 \\ \hline
\end{tabular}
\label{tab:example-D-Noisy-Solution-appendix}
\end{table}

\begin{table}[]
\caption{The experimental results for various methods under different metrics are presented for the FOP defined as $\min \{ x^{\top}x - (\theta + u)^{\top} x \mid x \in [0, 1]^{p} \} = \min \left\{\sum_{k=1}^{p} x_k^2 - \sum_{k=1}^{p} (\theta + u)_k x_k \mid x_k \in [0, 1], \ \forall \ k \in [p] \right\}$ in the \textbf{Noisy Objective Function} setting. For the FY loss, we set $\lambda=0.1$ and $\Omega(x) = 1/2\|x\|_2^2$.}
\vskip 0.15in
\centering
\begin{tabular}{c|cccc|cccc|cccc}
\hline
\multirow{2}{*}{Samplesize} & \multicolumn{4}{c|}{ParameterError} & \multicolumn{4}{c|}{DecisionError} & \multicolumn{4}{c}{Regret} \\ \cline{2-13} 
 & FY & SPA & KKA & VIA & FY & SPA & KKA & VIA & FY & SPA & KKA & VIA \\ \hline
50 & \textbf{1.37} & 2.52 & - & 8.80 & \textbf{0.05} & 0.15 & - & 0.85 & \textbf{0.05} & 0.16 & - & 0.85 \\
100 & \textbf{1.26} & 2.67 & - & 8.29 & \textbf{0.04} & 0.15 & - & 0.80 & \textbf{0.04} & 0.16 & - & 0.80 \\
300 & \textbf{1.28} & 2.88 & - & 5.85 & \textbf{0.03} & 0.17 & - & 0.51 & \textbf{0.04} & 0.18 & - & 0.51 \\
500 & \textbf{1.28} & 2.89 & - & 4.57 & \textbf{0.03} & 0.17 & - & 0.37 & \textbf{0.03} & 0.18 & - & 0.37 \\
1000 & \textbf{1.28} & 2.89 & - & 3.38 & \textbf{0.03} & 0.17 & - & 0.23 & \textbf{0.03} & 0.18 & - & 0.23 \\ \hline
\end{tabular}
\label{tab:example-D-Noisy-objective-function-appendix}
\end{table}

\begin{table}[]
\caption{The experimental results for various methods under different metrics are presented for the FOP defined as $\min \{ x^{\top}x - (\theta + u)^{\top} x \mid x \in [0, 1]^{p} \} = \min \left\{\sum_{k=1}^{p} x_k^2 - \sum_{k=1}^{p} (\theta + u)_k x_k \mid x_k \in [0, 1], \ \forall \ k \in [p] \right\}$ under \textbf{Noiseless} setting. For the FY loss, we set $\lambda=0.1$ and $\Omega(x) = 1/2 \|x\|_2^2$.}
\vskip 0.15in
\centering
\begin{tabular}{c|cccc|cccc|cccc}
\hline
\multirow{2}{*}{Sample size} & \multicolumn{4}{c|}{ParameterError} & \multicolumn{4}{c|}{DecisionError} & \multicolumn{4}{c}{Regret} \\ \cline{2-13} 
 & FY & SPA & KKA & VIA & FY & SPA & KKA & VIA & FY & SPA & KKA & VIA \\ \hline
50 & \textbf{0.07} & 1.26 & 0.00 & 13.29 & \textbf{0.00} & 0.03 & 0.00 & 1.33 & \textbf{0.00} & 0.03 & 0.00 & 1.33 \\
100 & \textbf{0.06} & 1.31 & 0.00 & 13.53 & \textbf{0.00} & 0.03 & 0.00 & 1.35 & \textbf{0.00} & 0.04 & 0.00 & 1.35 \\
300 & \textbf{0.06} & 1.32 & 0.00 & 14.66 & \textbf{0.00} & 0.03 & 0.00 & 1.40 & \textbf{0.00} & 0.04 & 0.00 & 1.40 \\
500 & \textbf{0.06} & 1.32 & 0.00 & 14.82 & \textbf{0.00} & 0.03 & 0.00 & 1.41 & \textbf{0.00} & 0.04 & 0.00 & 1.41 \\
1000 & \textbf{0.06} & 1.33 & 0.00 & 14.88 & \textbf{0.00} & 0.03 & 0.00 & 1.41 & \textbf{0.00} & 0.04 & 0.00 & 1.41 \\ \hline
\end{tabular}
\label{tab:example-D-Noiseless-appendix}
\end{table}

\clearpage
\subsubsection{Example E}
Finally, we consider the following forward optimization problem(FOP) with nonlinear constraint,
$$
    \min \{ -(\theta + u)^{\top} x \mid \parallel x \parallel_2^2 \leq a^2 \} = \min \left\{-\sum_{k=1}^{p} (\theta_k + u_k) x_k \mid \parallel x \parallel_2^2 \leq a^2 \right\}.
$$
where context vector $u$ follows a uniform distribution with support $[0,2]^{p}$ and $p=10$. We set the true parameter vector $\theta=[0.5,\cdots, 0.5] \in \mathbb{R}^p$ and $a=3$. The optimal solution to this problem can be expressed as
$$
    x^{\star}(\theta; u) = \frac{a(\theta+u)}{||\theta+u||_2}
$$
To compute the Fenchel-Young loss, we first derive the convex conjugate of $\Omega$, which is
$$
    \max_{x \in \mathcal{X}(u)} V_{\lambda}(\theta; u, x) = \max_{x \in \mathcal{F} } \left\{ (\theta + u)^{\top} x - \frac{\lambda}{2} x^{\top}x \right\} ~~~ \rightarrow ~~~ x^{\star}_{\lambda}(\theta; u) = \left\{
    \begin{aligned}
        & \frac{\theta+u}{\lambda}, && \text{if } \lambda \geq \frac{||\theta+u||_2}{a} \\
        & \frac{a(\theta+u)}{||\theta+u||_2}, && \text{if } \lambda < \frac{||\theta+u||_2}{a}.
    \end{aligned}
    \right.
$$
where feasible set $\mathcal{F}:\{x \in \mathbb{R}^p| \sum_{i=1}^p x_i^2 \leq a^2 \}$. Based on the observation solution $y$, we have the following expression of Fenchel-Young loss,
$$
\begin{aligned}
    L_{\lambda}(\theta; u, y) &= \max_{x \in \mathcal{X}(u)} V_{\lambda}(\theta; u, x) + \lambda\Omega(y) - \left[ \sum_{k=1}^{p} (\theta_k + u_k) y_{k} \right] \\ 
    &= \max_{x \in \mathcal{X}(u)} V_{\lambda}(\theta; u, x) + \lambda\Omega(y) - \sum_{k=1}^{p} (\theta_k + u_k) y_{k}
\end{aligned}
$$
Thus, we have the gradient of Fenchel-Young Loss:
$$
    \nabla_{\theta_k} L_{\lambda}(\theta;u, y) =  \big[x_{\lambda}^{\star}(\theta; u)_k - y_k \big]
$$
for $k = 1, \dots, p$. 
\paragraph{Result Discussion} The results of experiments in the Noisy Decision setting, the Noisy Objective Function setting, and the Noiseless setting, which are shown in \Cref{tab:example-E-Noisy-Solution-appendix}, \Cref{tab:example-E-Noisy-objective-function-appendix} and \Cref{tab:example-E-Noiseless-appendix} respectively, average over 100 repetitions. Based on this information, we can draw the following conclusions.
\begin{itemize}
    \item In the Noisy Decision and Noisy Objective Function settings, only our method based on the Fenchel-Young loss can be effectively solved using Gurobi, while other inverse optimization methods, except KKA, fail to provide effective solutions. Although KKA can be solved, the results do not converge. This highlights the robustness and efficiency of our Fenchel-Young loss method in handling inverse optimization problems with noisy data and nonlinear constraints, which are common challenges in this domain.

    \item In the Noiseless setting, the VIA method demonstrates the best performance, achieving a perfect estimation of the parameter vector $\theta$ when the sample size was equal to 50. Furthermore, both its decision error and regret are zero. However, as the sample size increases, the computational time exceeds practical limits, preventing the acquisition of results. Consequently, these cases were not recorded in the table.
    
    In addition, KKA fails to converge under the Noiseless setting, resulting in substantial parameter error, decision error, and regret. However, our proposed method continues to deliver near-perfect results with rapid convergence. Specifically, the parameter error is remarkably low, while the decision error and regret are both zero.
\end{itemize}

\begin{table}[ht]
\caption{The experimental results for various methods under different metrics are presented for the FOP defined as $\min \{ -(\theta + u)^{\top} x \mid \parallel x \parallel_2^2 \leq a^2 \} = \min \left\{-\sum_{k=1}^{p} (\theta_k + u_k) x_k \mid \parallel x \parallel_2^2 \leq a^2 \right\}$ in hte \textbf{Noisy Decision} setting. For the FY loss, we set $\lambda=0.1$ and $\Omega(x) = 1/2 \|x\|_2^2$.}
\vskip 0.15in
\centering
\begin{tabular}{c|cccc|cccc|cccc}
\hline
\multirow{2}{*}{Sample size} & \multicolumn{4}{c|}{Parameter Error} & \multicolumn{4}{c|}{Decision Error} & \multicolumn{4}{c}{Regret} \\ \cline{2-13} 
 & FY & SPA & KKA & VIA & FY & SPA & KKA & VIA & FY & SPA & KKA & VIA \\ \hline
50 & \textbf{1.02} & - & 4.59 & - & \textbf{0.24} & - & 1.08 & - & \textbf{0.09} & - & 0.38 & - \\
100 & \textbf{0.67} & - & 4.80 & - & \textbf{0.11} & - & 1.13 & - & \textbf{0.04} & - & 0.40 & - \\
300 & \textbf{0.39} & - & 4.93 & - & \textbf{0.04} & - & 1.16 & - & \textbf{0.01} & - & 0.41 & - \\
500 & \textbf{0.30} & - & 4.96 & - & \textbf{0.02} & - & 1.17 & - & \textbf{0.01} & - & 0.42 & - \\
1000 & \textbf{0.23} & - & 4.98 & - & \textbf{0.01} & - & 1.18 & - & \textbf{0.00} & - & 0.42 & - \\ \hline
\end{tabular}
\label{tab:example-E-Noisy-Solution-appendix}
\end{table}

\begin{table}[ht]
\caption{The experimental results for various methods under different metrics are presented for the FOP defined as $\min \{ -(\theta + u)^{\top} x \mid \parallel x \parallel_2^2 \leq a^2 \} = \min \left\{-\sum_{k=1}^{p} (\theta_k + u_k) x_k \mid \parallel x \parallel_2^2 \leq a^2 \right\}$ in the \textbf{Noisy Objective Function} setting. For the FY loss, we set $\lambda=0.1$ and $\Omega(x) = 1/2\|x\|_2^2$.}
\vskip 0.15in
\centering
\begin{tabular}{c|cccc|cccc|cccc}
\hline
\multirow{2}{*}{Sample size} & \multicolumn{4}{c|}{Parameter Error} & \multicolumn{4}{c|}{Decision Error} & \multicolumn{4}{c}{Regret} \\ \cline{2-13} 
 & FY & SPA & KKA & VIA & FY & SPA & KKA & VIA & FY & SPA & KKA & VIA \\ \hline
50 & \textbf{0.51} & - & 4.59 & - & \textbf{0.04} & - & 1.08 & - & \textbf{0.01} & - & 0.38 & - \\
100 & \textbf{0.50} & - & 4.80 & - & \textbf{0.03} & - & 1.13 & - & \textbf{0.01} & - & 0.40 & - \\
300 & \textbf{0.49} & - & 4.93 & - & \textbf{0.03} & - & 1.16 & - & \textbf{0.01} & - & 0.41 & - \\
500 & \textbf{0.50} & - & 4.96 & - & \textbf{0.03} & - & 1.17 & - & \textbf{0.01} & - & 0.42 & - \\
1000 & \textbf{0.50} & - & 4.98 & - & \textbf{0.03} & - & 1.18 & - & \textbf{0.01} & - & 0.42 & - \\ \hline
\end{tabular}
\label{tab:example-E-Noisy-objective-function-appendix}
\end{table}

\begin{table}[ht]
\caption{The experimental results for various methods under different metrics are presented for the FOP defined as $\min \{ -(\theta + u)^{\top} x \mid \parallel x \parallel_2^2 \leq a^2 \} = \min \left\{-\sum_{k=1}^{p} (\theta_k + u_k) x_k \mid \parallel x \parallel_2^2 \leq a^2 \right\}$ in the \textbf{Noiseless} setting. For the FY loss, we set $\lambda=0.1$ and $\Omega(x) = 1/2\|x\|_2^2$}
\vskip 0.15in
\centering
\begin{tabular}{c|cccc|cccc|cccc}
\hline
\multirow{2}{*}{Sample size} & \multicolumn{4}{c|}{Parameter Error} & \multicolumn{4}{c|}{Decision Error} & \multicolumn{4}{c}{Regret} \\ \cline{2-13} 
 & FY & SPA & KKA & VIA & FY & SPA & KKA & VIA & FY & SPA & KKA & VIA \\ \hline
50 & \textbf{0.04} & - & 4.59 & 0.00 & \textbf{0.00} & - & 1.08 & 0.00 & \textbf{0.00} & - & 0.38 & 0.00 \\
100 & \textbf{0.04} & - & 4.80 & - & \textbf{0.00} & - & 1.13 & - & \textbf{0.00} & - & 0.40 & - \\
300 & \textbf{0.04} & - & 4.93 & - & \textbf{0.00} & - & 1.16 & - & \textbf{0.00} & - & 0.41 & - \\
500 & \textbf{0.04} & - & 4.96 & - & \textbf{0.00} & - & 1.17 & - & \textbf{0.00} & - & 0.42 & - \\
1000 & \textbf{0.03} & - & 4.98 & - & \textbf{0.00} & - & 1.18 & - & \textbf{0.00} & - & 0.42 & - \\ \hline
\end{tabular}
\label{tab:example-E-Noiseless-appendix}
\end{table}
\clearpage

\subsubsection{Computational Time} \label{sec:computational-time-appendix}
In this section, we present a concise summary of the total computational time for 100 serial runs of the synthetic data experiments under various settings and sample sizes. Taking the Noisy Decision setting as an example, the results for other settings are similar, and the detailed findings are shown in \Cref{tab:Computing-time-results-Noisy-solution-appendix}. The units of the numbers in the table are seconds. From these results, the following conclusions can be readily drawn.
\begin{itemize}
    \item \textbf{Computational Efficiency of Our Method:} Except for Example A, our method consistently achieves the shortest computational time across all experiments with an order-of-magnitude improvement, typically requiring only a few seconds to complete the entire computation in most experiments. In contrast, other methods generally take at least tens of seconds or more. Additionally, as the sample size increases, the growth rate of computational time for our method remains manageable. In Example A, the explicit gradient expression of the Fenchel-Young loss cannot be derived, necessitating the use of a solver, which slows down the computation. However, by tuning the parameters of stochastic gradient descent, we can still complete the computation within an acceptable time frame, even with a large sample size.
    
    \item \textbf{Computational Bottleneck of SPA:} In all experiments, SPA exhibits the longest computational time among all methods. This is attributed to its use of a semi-parametric algorithm for noise removal. According to relevant literature \cite{aswani2018inverse}, the semi-parametric algorithm requires cross-validation for hyperparameter tuning, which significantly increases the computational time compared to KKA and VIA.

    \item \textbf{Comparable Computational Time of KKA and VIA: } The computational time for KKA and VIA remains on the same order of magnitude across all experiments. Except for Example A, their computational times are significantly higher than that of our method but considerably lower than that of SPA.
\end{itemize}

\begin{table}[]
\caption{computational time results of different methods in each experiments under \textbf{Noisy Decision} setting.}
\vskip 0.15in
\centering
\begin{tabular}{cc|ccccc}
\hline
\multicolumn{2}{c|}{Sample size} & \multicolumn{1}{c|}{50} & \multicolumn{1}{c|}{100} & \multicolumn{1}{c|}{300} & \multicolumn{1}{c|}{500} & 1000 \\ \hline
\multicolumn{1}{c|}{\multirow{4}{*}{Example A}} & FY & \textbf{20.34} & \textbf{21.54} & \textbf{68.33} & \textbf{64.13} & \textbf{97.20} \\ \cline{2-2}
\multicolumn{1}{c|}{} & SPA & 290.84 & 409.33 & 710.18 & 1858.89 & 3574.48 \\ \cline{2-2}
\multicolumn{1}{c|}{} & KKA & 2.92 & 9.76 & 13.76 & 26.71 & 49.82 \\ \cline{2-2}
\multicolumn{1}{c|}{} & VIA & 2.19 & 15.09 & 26.11 & 41.70 & 78.80 \\ \hline
\multicolumn{1}{c|}{\multirow{4}{*}{Example B}} & FY & \textbf{0.75} & \textbf{0.85} & \textbf{1.54} & \textbf{2.28} & \textbf{4.39} \\ \cline{2-2}
\multicolumn{1}{c|}{} & SPA & 62.54 & 122.73 & 504.12 & 1064.74 & 3906.52 \\ \cline{2-2}
\multicolumn{1}{c|}{} & KKA & 3.09 & 6.52 & 20.73 & 37.43 & 70.86 \\ \cline{2-2}
\multicolumn{1}{c|}{} & VIA & 2.41 & 3.42 & 32.94 & 37.02 & 36.83 \\ \hline
\multicolumn{1}{c|}{\multirow{4}{*}{Example C}} & FY & \textbf{0.72} & \textbf{0.77} & \textbf{1.25} & \textbf{1.85} & \textbf{3.25} \\ \cline{2-2}
\multicolumn{1}{c|}{} & SPA & 166.96 & 140.48 & 507.64 & 1024.51 & 3205.72 \\ \cline{2-2}
\multicolumn{1}{c|}{} & KKA & 3.38 & 7.78 & 20.12 & 32.53 & 62.99 \\ \cline{2-2}
\multicolumn{1}{c|}{} & VIA & 1.98 & 8.02 & 15.05 & 24.34 & 54.06 \\ \hline
\multicolumn{1}{c|}{\multirow{4}{*}{Example D}} & FY & \textbf{0.65} & \textbf{0.76} & \textbf{1.30} & \textbf{1.97} & \textbf{3.43} \\ \cline{2-2}
\multicolumn{1}{c|}{} & SPA & 86.24 & 159.85 & 590.19 & 1780.07 & 4481.71 \\ \cline{2-2}
\multicolumn{1}{c|}{} & KKA & - & - & - & - & - \\ \cline{2-2}
\multicolumn{1}{c|}{} & VIA & 13.62 & 9.93 & 38.16 & 76.80 & 204.85 \\ \hline
\multicolumn{1}{c|}{\multirow{4}{*}{Example E}} & FY & \textbf{1.94} & \textbf{3.73} & \textbf{14.98} & \textbf{35.67} & \textbf{38.22} \\ \cline{2-2}
\multicolumn{1}{c|}{} & SPA & - & - & - & - & - \\ \cline{2-2}
\multicolumn{1}{c|}{} & KKA & 42.46 & 257.40 & 907.82 & 2486.15 & 21449.75 \\ \cline{2-2}
\multicolumn{1}{c|}{} & VIA & - & - & - & - & - \\ \hline
\end{tabular}
\label{tab:Computing-time-results-Noisy-solution-appendix}
\end{table}
\clearpage

\begin{figure}[ht]
    \centering
    \begin{minipage}[b]{0.325\linewidth}
        \centering
        \includegraphics[width=1\linewidth]{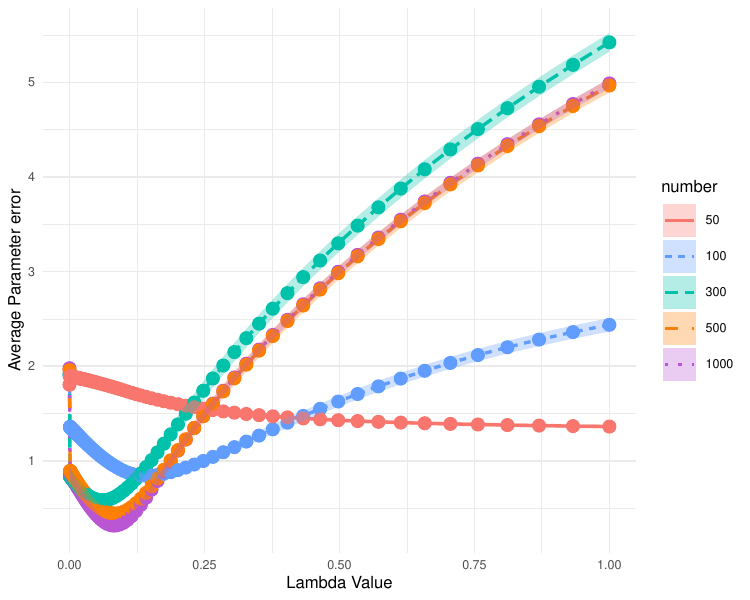}
        \caption{Example A PE}
        \label{fig:Example A parameter error}
    \end{minipage}
    \hfill
    \begin{minipage}[b]{0.325\linewidth}
        \centering
        \includegraphics[width=1\linewidth]{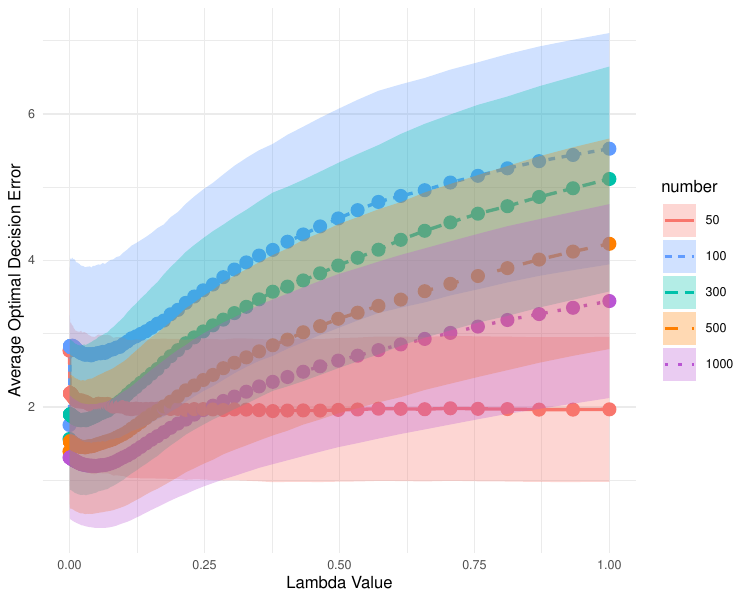}
        \caption{Example A DE}
    \end{minipage}
    \hfill
    \begin{minipage}[b]{0.325\linewidth}
        \centering
        \includegraphics[width=1\linewidth]{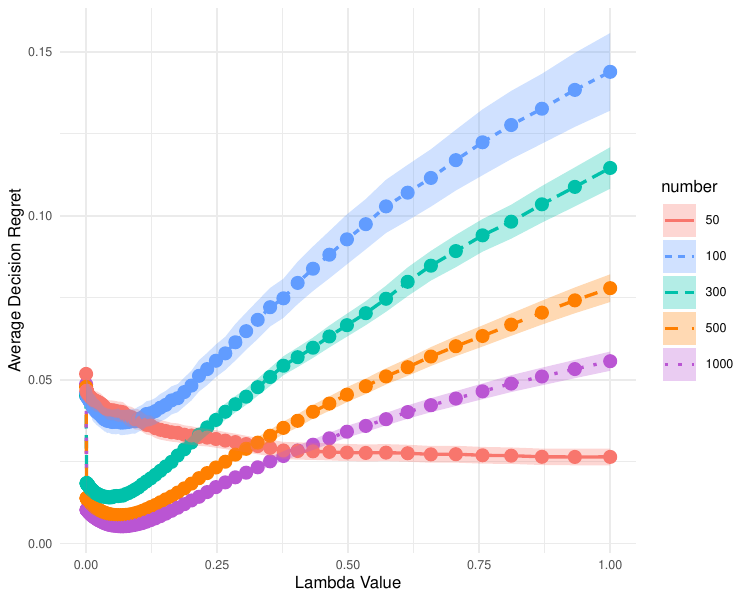}
        \caption{Example A regret}
    \end{minipage}
    \qquad
    \centering
    \begin{minipage}[b]{0.325\linewidth}
        \centering
        \includegraphics[width=1\linewidth]{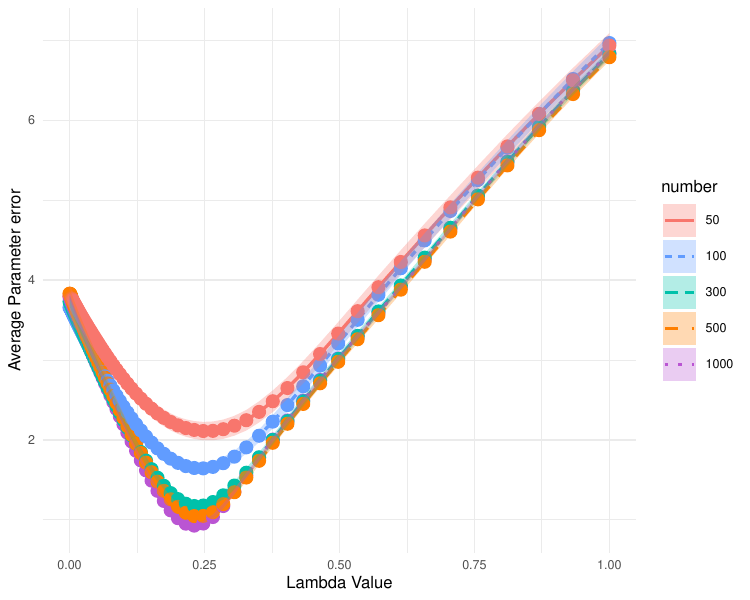}
        \caption{Example B PE}
        \label{fig:Example B parameter error}
    \end{minipage}
    \hfill
    \begin{minipage}[b]{0.325\linewidth}
        \centering
        \includegraphics[width=1\linewidth]{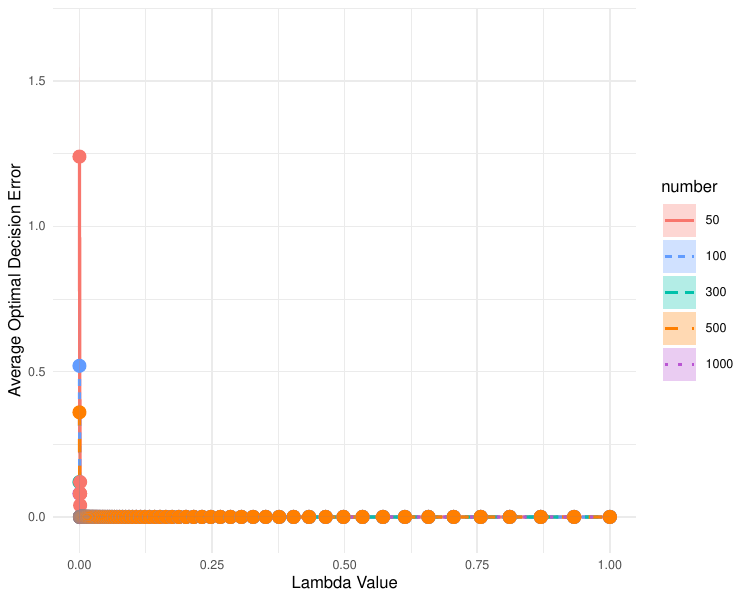}
        \caption{Example B DE}
    \end{minipage}
    \hfill
    \begin{minipage}[b]{0.325\linewidth}
        \centering
        \includegraphics[width=1\linewidth]{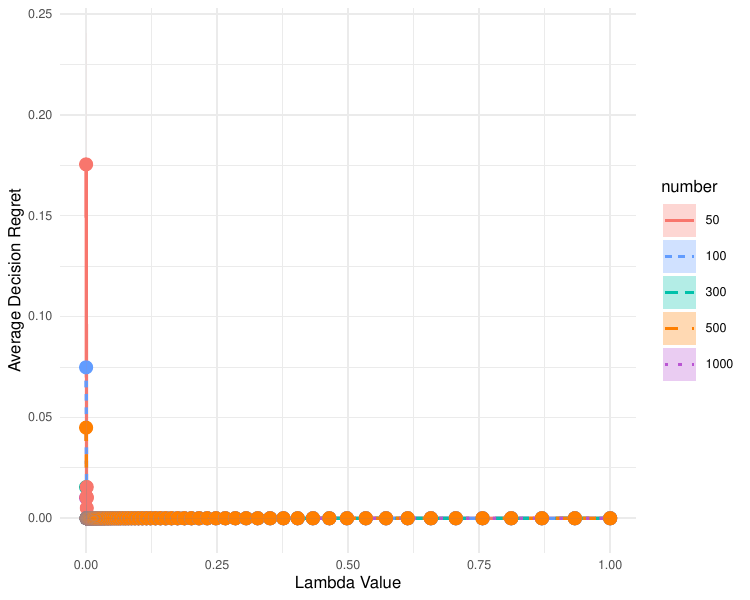}
        \caption{Example B regret}
    \end{minipage}
    \qquad

    \centering
    \begin{minipage}[b]{0.325\linewidth}
        \centering
        \includegraphics[width=1\linewidth]{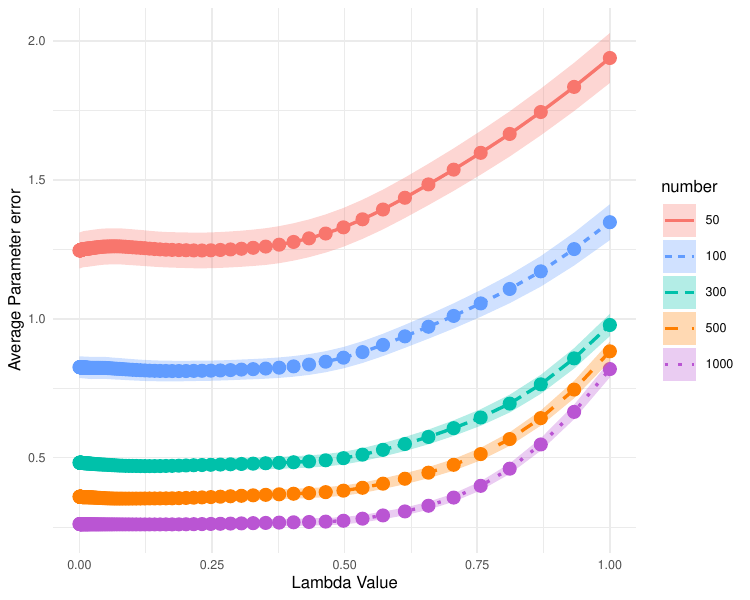}
        \caption{Example C PE}
    \end{minipage}
    \hfill
    \begin{minipage}[b]{0.325\linewidth}
        \centering
        \includegraphics[width=1\linewidth]{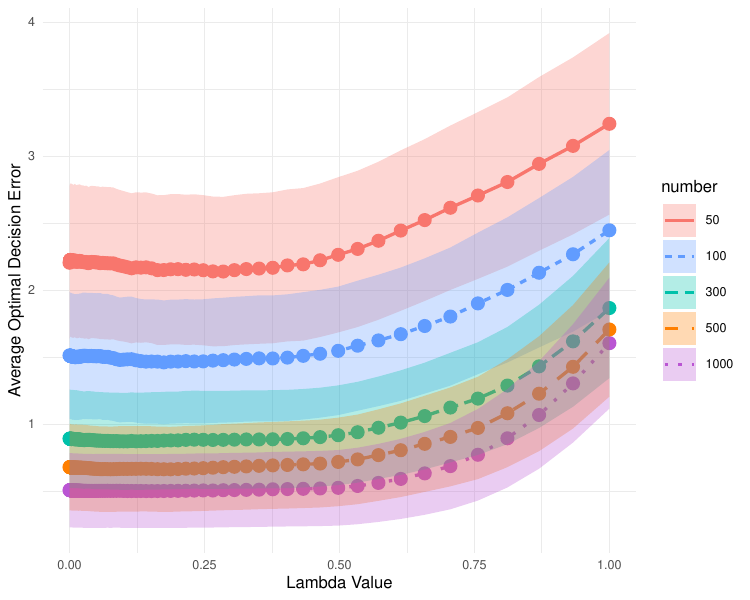}
        \caption{Example C DE}
    \end{minipage}
    \hfill
    \begin{minipage}[b]{0.33\linewidth}
        \centering 
        \includegraphics[width=1\linewidth]{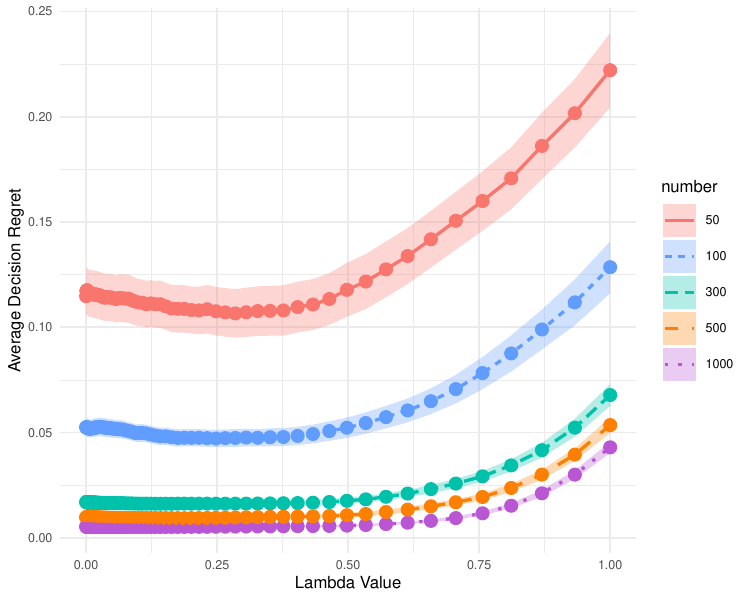}
        \caption{Example C regret}
    \end{minipage}
\end{figure}

\begin{figure}[ht]
    \centering
    \begin{minipage}[b]{0.325\linewidth}
        \centering
        \includegraphics[width=1\linewidth]{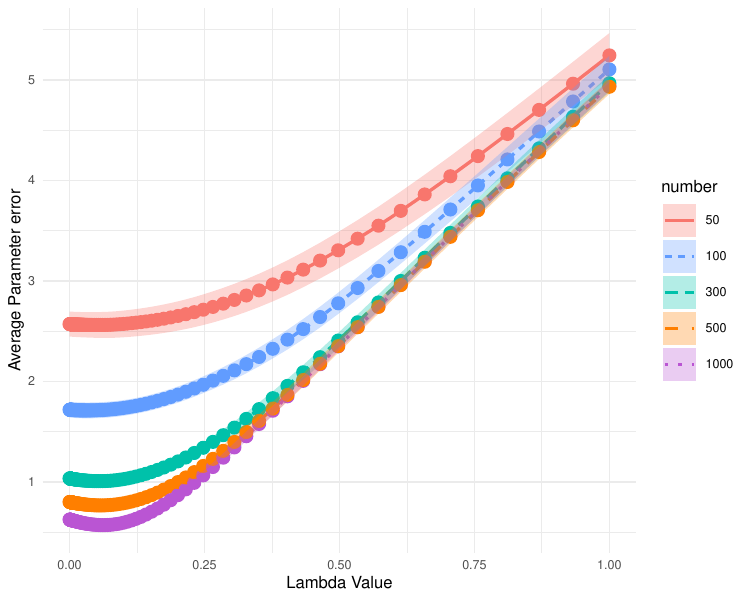}
        \caption{Example D PE}
    \end{minipage}
    \hfill
    \begin{minipage}[b]{0.325\linewidth}
        \centering
        \includegraphics[width=1\linewidth]{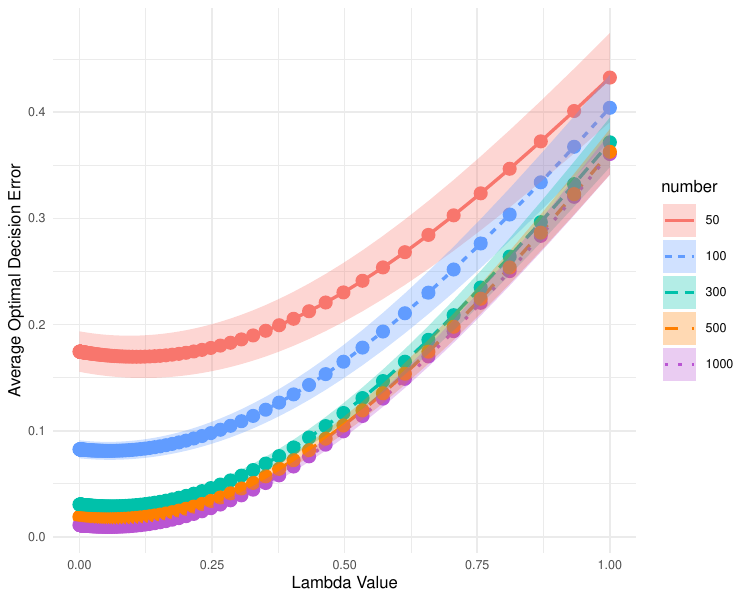}
        \caption{Example D DE}
    \end{minipage}
    \hfill
    \begin{minipage}[b]{0.325\linewidth}
        \centering
        \includegraphics[width=1\linewidth]{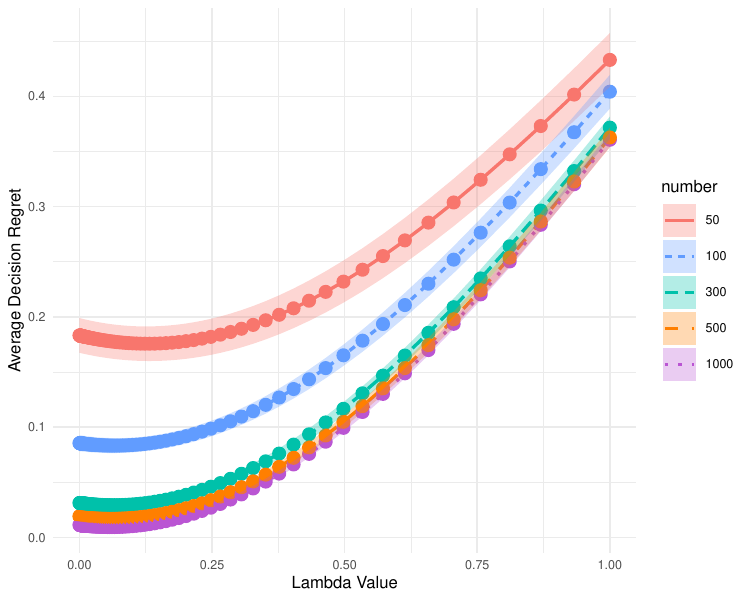}
        \caption{Example D regret}
    \end{minipage}
    \qquad
    \centering
    \begin{minipage}[b]{0.325\linewidth}
        \centering
        \includegraphics[width=1\linewidth]{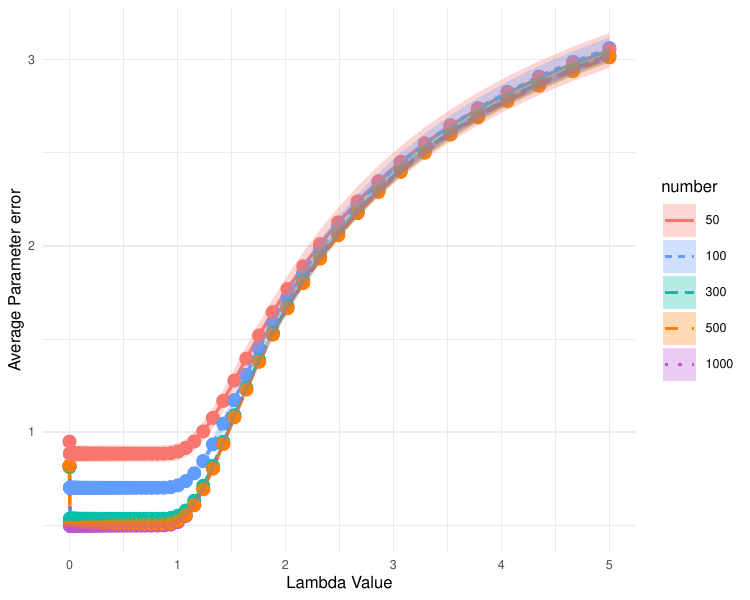}
        \caption{Example E PE}
    \end{minipage}
    \hfill
    \begin{minipage}[b]{0.325\linewidth}
        \centering
        \includegraphics[width=1\linewidth]{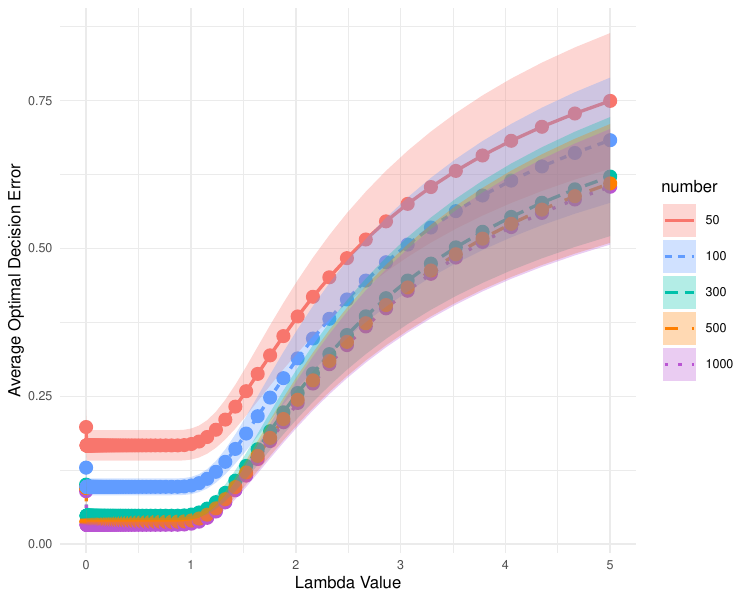}
        \caption{Example E DE}
    \end{minipage}
    \hfill
    \begin{minipage}[b]{0.325\linewidth}
        \centering
        \includegraphics[width=1\linewidth]{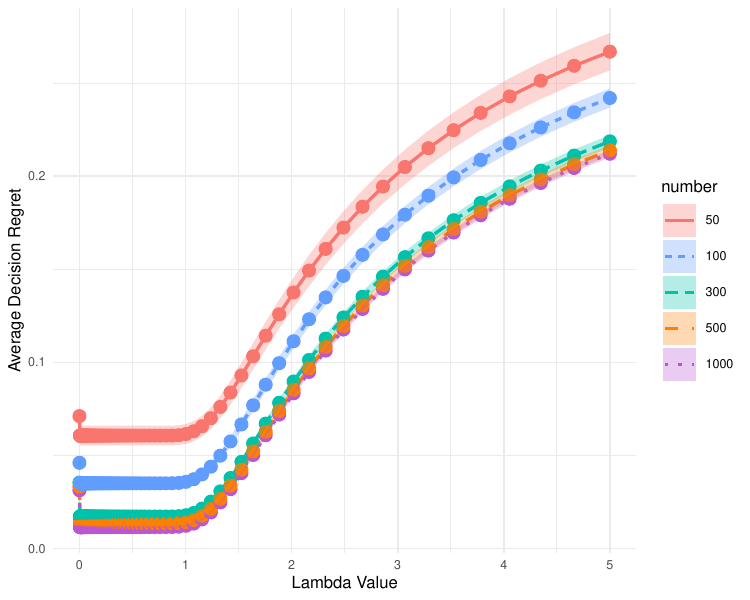}
        \caption{Example E regret}
        \label{fig:Example E regret}
    \end{minipage}
\end{figure}

\subsubsection{Sensitivity Analysis}
\label{sec:sensitivity-analysis-appendix}
In this section, we investigate the impact of the regularization parameter $\lambda$ on various metrics when solving inverse optimization problems using our method. Specifically, for all experiments except Example E, we select regular coefficient $\lambda$ from a grid of 0 and 100 points distributed uniformly on the logarithmic scale over 0.001 to 1. For Example E, we select regular coefficient $\lambda$ from a grid of 0 and 100 points distributed uniformly on the logarithmic scale over 0.005 to 5. 

From \Cref{fig:Example A parameter error} to \Cref{fig:Example E regret}, we present the results of the Fenchel-Young loss for different $\lambda$ values in each synthetic data experiment under the Noisy Decision setting, along with the corresponding 95\% confidence intervals derived from 100 repeated experiments. In the titles of these figures, \textbf{PE} denotes Parameter Error, while \textbf{DE} denotes Decision Error. Based on these results, we draw the following conclusions.
\begin{itemize}
    \item For Example A, When the sample size is 50, the parameter error, decision error, and regret gradually decrease as $\lambda$ increases. This indicates that when the sample size is insufficient, a sufficiently strong regularization is beneficial for better performance. When the sample size exceeds 50, these three metrics first decrease and then increase with the increase of $\lambda$, showing a V-shaped pattern. This suggests that when the sample size is adequate, excessive regularization introduces significant bias, which is detrimental to performance enhancement. For parameter error, the optimal $\lambda$ value varies with different sample sizes. In contrast, for decision error and regret, the optimal $\lambda$ value remains almost the same across different sample sizes. 
    
    In summary, this example primarily illustrates that the selection of the optimal $\lambda$ value requires a comprehensive consideration of multiple factors, including sample size and specific metrics. Moreover, in all three figures, there is a sharp drop when $\lambda$ transitions from 0 to a non-zero value, highlighting the necessity and rationality of regularization.

    \item For Example B, the parameter error curve exhibits a distinct V-shape. Specifically, as the value of $\lambda$ increases, the parameter error initially decreases and then subsequently increases. The minimum point is approximately at $\lambda=0.25$, and as the sample size increases, this minimum point shifts slightly to the left. This trend is reasonable because, with a larger sample size, the optimal value of the regularization term should gradually decrease.

    Moreover, in this example, the decision error and regret curves clearly illustrate the benefits of regularization. When $\lambda=0$, the solution obtained is far from optimal, and the regret is substantial. However, once $\lambda>0$, the solution immediately becomes nearly optimal.

    In summary, this example demonstrates that regularization is both helpful and necessary. Even a small value of $\lambda$ can significantly outperform a value of 0. Additionally, proper selection of $\lambda$ can facilitate parameter recovery.

    \item For Example C, it is observed that when the value of $\lambda$ lies between 0 and 0.5, the three metrics under consideration remain relatively stable overall. Specifically, when the sample size is 50, the regret exhibits a slight trend of initially decreasing and then increasing, with the lowest point occurring around $\lambda = 0.3$.
    
    However, when the value of $\lambda$ exceeds 0.5, the three metrics show a significant upward trend. This is reasonable because an excessively large regularization coefficient $\lambda$ introduces substantial bias, leading to poor performance. Additionally, as the sample size increases, the values of the three metrics decrease, which is expected.

    In summary, this example demonstrates that our method allows for a broad range of $\lambda$ values. As long as $\lambda$ is not excessively large, it can ensure effective parameter recovery and satisfactory decision-making performance.
    
    \item For Example D, compared to Example C, its objective function includes a quadratic term $x^{\top}x$, and the term $(\theta+u)$ is modified to $-(\theta+u)$. Similar to Example C, when the $\lambda$ value is very small, specifically around 0 to 0.1, the curves for parameter error, decision error, and regret are relatively smooth. However, when the value of $\lambda$ exceeds 0.1, these metrics increase rapidly. Additionally, when the sample size is 50, the regret initially decreases slightly before increasing. This may be because, when the sample size is insufficient, appropriate regularization helps to reduce regret.
    
    In summary, this example illustrates that when dealing with nonlinear programming, the choice of $\lambda$ may need to be more carefully considered compared to linear programming, with values closer to 0 being more appropriate.
 
    \item For Example E, we specifically adjust the $\lambda$ value selection range to 0.005 to 5 instead of 0.001 to 1, just to show the two phase transitions in this experiment. First, we find that there is a jump at the point where $\lambda = 0$. In other words, when $\lambda = 0$ becomes $\lambda > 0$, parameter error, decision error and regret will suddenly decrease, which reflects the advantages of regularization. For $\lambda$ values are less than 1, the three metrics we consider remain at a very stable and very low level. Reviewing the previous discussion of this example, we can see that the reason for this is that when $\lambda$ is less than $\frac{||\theta+u||_2}{a}$, $x^*(\theta,u)=\frac{a(\theta+u)}{||\theta+u||_2}$ is a fixed constant and our algorithm is not affected in any way. Only when $\lambda$ is greater than $\frac{||\theta+u||_2}{a}$ does $x^*(\theta,u)=\frac{\theta+u}{\lambda}$ change with $\lambda$. Consequently, when $\lambda$ is greater than 1, the phase transition occurs again. Due to the introduction of large bias, parameter error, decision error and regret will rise rapidly with the increase of $\lambda$.

    In general, this example demonstrates the unique properties of our method in the ball constraint space. When regularization is present and the value is below a certain threshold, our method ensures a very stable and relatively optimal performance.

\end{itemize}

\subsubsection{Synthetic Data Experiment Conclusions}

In summary, our method is gradient-based, whereas other classical inverse optimization methods are primarily model-based. This fundamental difference endows our approach with unique advantages in several key aspects:
\begin{itemize}
    \item \textbf{Avoid Degeneracy}: When optimization parameters are estimated using model-based inverse optimization methods, degenerate solutions (e.g., an all-zero vector in Example B) may arise, leading to meaningless optimal solutions and ineffective decision-making. Our Fenchel-Young loss-based method circumvents this issue, ensuring robust decision performance.
    \item \textbf{Handling Noise More Effectively}: Many classical inverse optimization methods perform well only in noiseless settings. In noisy environments, solver failures are common. In contrast, our method remains stable regardless of noise, yielding accurate parameter estimates and satisfactory decision outcomes.
    \item\textbf{Effective Acceleration}: Our method exhibits significant computational advantages over other inverse optimization methods. Although we cannot derive an explicit expression for the optimal solution $x^{\star}_{\lambda}(\theta;u)$ of the conjugate function $\max_{x \in \mathcal{X}(u)} V_{\lambda}(\theta; u, x)$ in Example A, we can employ stochastic gradient descent to maintain acceptable computational time, even with large sample sizes. Conversely, the computational time for other classical methods increases sharply with increasing sample size. These methods, being model-based, incorporate each data point as a constraint, thereby exponentially increasing the complexity of optimization problem. Moreover, there is no systematic approach to accelerate these calculations.
\end{itemize}
\clearpage

\subsection{Real Data Experiment Setup and Implementation details} \label{sub:real-data-experiment}
In this section we study our Fenchel-Young loss method and baselines empirically to demonstrate the great advantages of our method in solving inverse optimization problems in practical applications. The real data we use are collected from Uber Movement(https://movement.uber.com), which describes the traveling times data in Los Angeles. In our experiment, We mainly focus on 45 census tracts in downtown Los Angeles, collecting historical data of average traveling times from each of these census tracts to its neighbors, which includes 93 edges, during five periods in each day (AM Peak, Midday, PM Peak, Evening, Early Morning) in 2018 and 2019. Next, we will introduce our experiment from the following aspects.
\begin{itemize}
    \item \textbf{Forward optimization problem(FOP) and observed solution}: Our objective is to determine an optimal path from the easternmost census tract(Aliso Village) to the westernmost census tract(MacArthur Park), which is encoded by $x \in \{0,1\}^d$ with $d=93$. More detailed setting can be found in \cite{kallus2023stochastic}. When $x_j=1$, it means that we travel on edge $j$. Because our forward optimization problem is actually the shortest path problem on a graph with 45 nodes and 93 edges. Therefore, we also need to consider the flow preservation constraints, which are denoted as $Ax=b$,$A \in \mathbb{R}^{45 \times 93}$,$b \in \mathbb{R}^{45}$.
    
    Furthermore, to construct an inverse optimization problem, we solve for a shortest path based on the travel time of each edge at each data point in the dataset, which serves as the observed optimal solution $Y$.  The travel times vector for 93 edges is modeled as $\boldsymbol{\theta} u$, where $\boldsymbol{\theta} \in \mathbb{R}^{93 \times 12}$ and $u \in \mathbb{R}^{12}$. 12 means that the dimension of the context vector of each piece of data is 12, which we will describe in detail in the next paragraph. Finally, the forward optimization problem(FOP) is shown below.

    $$
    \min_{x} \left\{\sum_{k=1}^{93} (\boldsymbol{\theta}_k u) x_k \mid Ax=b \right\}.
    $$

    $\boldsymbol{\theta}_k$ is the k-th row of the matrix $\boldsymbol{\theta}$. The objective of inverse optimization is to find a parameter matrix $\boldsymbol{\theta} \in \mathbb{R}^{93 \times 12}$ and estimate $c$ by $\boldsymbol{\theta} u$, so as to obtain a good decision on the test set data that is close enough to the observed optimal solution or low enough to the regret.
    
    \item \textbf{Context vector}: We totally consider 7 covariates including wind speed, visibility and other calendar features, which can be regarded as context vector $u$ for over inverse optimization problem. To be more precise, we process the period feature with one-hot encoder and add the intercept term, and the final processed context vector is $u \in \mathbb{R}^{12}$. 

    \item \textbf{Other experimental details}: We divide the datasets into four distinct time spans: half a year (920 data points), one year (1825 data points), one and a half years (2735 data points), and two years (3640 data points). For each dataset, we randomly select 40\% of the data as the test set and 60\% as the training set. Each set of experiments is run in parallel across 30 replications, and we collect the average computational time, the average relative regret ratio, and the average decision error.

    The relative regret ratio means that the regret on the test set is divided by the traveling time of the optimal decision, which can represent the percentage of time that the outcome of the decision takes longer than the optimal decision. The definition of decision error is exactly the same as in the previous synthetic data experiments.
\end{itemize}

Next, we delve into the specifics of how our method addresses these inverse optimization problems. 
To compute the Fenchel-Young loss, we first derive the convex conjugate of $\Omega$, which is
$$
    \max_{x \in \mathcal{X}(u)} V_{\lambda}(\boldsymbol{\theta}; u, x) = \max_{x \in \mathcal{F} } \left\{\sum_{k=1}^{93} (\boldsymbol{\theta}_k u) x_k - \frac{\lambda}{2} x^{\top}x \right\} ~~~ \rightarrow ~~~ x^{\star}_{\lambda}(\boldsymbol{\theta}; u).
$$
where feasible set $\mathcal{F}:\{x \in \{0,1\}^{93}| Ax=b \}$.  As in Example A, we still cannot directly write an explicit expression for x in this problem, so we need to utilize the solver to solve it. In order to ensure the efficiency of the gradient calculation, we still adopt the stochastic gradient descent algorithm, but different from Example A, we need to calculate that the gradient of $\boldsymbol{\theta}$ is a matrix rather than a vector, so we need to pay attention to the dimension matching. The details of gradient calculation are shown below.

Based on the observation solution $y$, we have the following expression of Fenchel-Young loss,
$$
\begin{aligned}
    L_{\lambda}(\boldsymbol{\theta};u,y) &= \max_{x \in \mathcal{X}(u)} V_{\lambda}(\boldsymbol{\theta}; u, x) + \lambda\Omega(y) - \left[ - \sum_{k=1}^{93} (\boldsymbol{\theta}_k u) y_k  \right] \\ 
    &= \max_{x \in \mathcal{X}(u)} V_{\lambda}(\boldsymbol{\theta}; u, x) + \lambda\Omega(y) + \sum_{k=1}^{93} (\boldsymbol{\theta}_k u) y_k
\end{aligned}
$$
Thus, we have the gradient of Fenchel-Young Loss:
$$
 \nabla_{\boldsymbol{\theta}} L_{\lambda}(\theta;u,y) =  \big[- x^{\star}_{\lambda}(\boldsymbol{\theta}; u) + y \big]u^{\top}
$$

\paragraph{Result Discussion} The results of this real data experiment are shown in \Cref{tab:real-data-experiment-result-appendix}. Occasionally, the VIA method fails to solve the problem. In such cases, we set $\boldsymbol{\theta}$ to be a matrix of all zeros, resulting in a random selection of a feasible path as the decision. Given that the forward optimization problem (FOP) for this problem is a linear programming problem and the observed solutions are noise-free, the maximum optimality margin (MOM) method is applicable. Based on the experimental results, we draw the following conclusions.
\begin{itemize}
    \item For decision error, our method achieves the lowest decision error, significantly outperforming other methods. This indicates that, regardless of the time span, the optimal path identified by our method closely approximates the true optimal path. The second-best performance is observed with the MOM method, which consistently ensures that the selected path is near-optimal. In contrast, the SPA, VIA, and KKA methods exhibit suboptimal performance, with substantial deviations from the true optimal solution.
    \item For relative regret ratio, our method exhibits significant advantages over other methods again. This indicates that the path travel time provided by our method is very close to the optimal path travel time. The second-best performance is observed with the MOM method, which also performs well and maintains a relatively low regret ratio. Although KKA has a high decision error, its relative regret ratio is marginally satisfactory. In contrast, the performance of SPA and VIA remains poor, with high relative regret ratios.
    \item For computational time, our method exhibits the shortest computational time compared to other methods, which are significantly longer. Other methods involve solving complex optimization problems with constraints that grow with the amount of data. Given the large volume of real data, these methods require extensive computational resources and time to reach a solution. In contrast, our method only requires a finite number of gradient descent steps, and computing the gradient involves solving relatively simple optimization problems. As a result, our method offers a substantial advantage in terms of computational speed.
\end{itemize}

\begin{table}[]
\caption{Real data experiment results. Numbers in the \textbf{Relative Regret Ratio} column mean percentages. For example, 1.01 means a relative regret ratio of 1.01\%. The number in the \textbf{Period} column means years. For example, 1.5 means that the data set used spans one and a half years. The number in the \textbf{Computational time} column means seconds. For example, 6 means that the average running time of the method is 6 seconds. For Fenchel-Young loss, we set $\lambda=0.1$ and $\Omega(x) = 1/2\|x\|_2^2$.}
\vskip 0.15in
\centering
\resizebox{\textwidth}{12mm}{
\begin{tabular}{c|clccc|ccccc|ccccc}
\hline
\multirow{2}{*}{\begin{tabular}[c]{@{}c@{}}Period\\ (years)\end{tabular}} & \multicolumn{5}{c|}{Decision Error} & \multicolumn{5}{c|}{Relative Regret Ratio (\%)} & \multicolumn{5}{c}{Computing time (seconds)} \\ \cline{2-16} 
 & FY & \multicolumn{1}{c}{SPA} & KKA & VIA & MOM & FY & SPA & KKA & VIA & MOM & FY & SP & KKA & VIA & MOM \\ \hline
0.5 & \textbf{2.06} & 6.80 & 5.41 & 10.72 & 3.58 & \textbf{1.01} & 29.51 & 8.30 & 34.51 & 4.87 & \textbf{6} & 120 & 66 & 59 & 438 \\
1 & \textbf{2.00} & 6.72 & 6.61 & 11.77 & 4.09 & \textbf{0.97} & 29.55 & 12.39 & 41.51 & 4.96 & \textbf{8} & 287 & 152 & 107 & 647 \\
1.5 & \textbf{1.97} & 6.77 & 7.09 & 7.33 & 4.66 & \textbf{0.94} & 29.50 & 12.81 & 30.95 & 4.21 & \textbf{9} & 513 & 242 & 185 & 901 \\
2 & \textbf{1.89} & 6.78 & 7.07 & 7.17 & 4.64 & \textbf{0.84} & 29.28 & 12.40 & 30.41 & 3.26 & \textbf{11} & 839 & 281 & 256 & 1218 \\ \hline
\end{tabular}}
\label{tab:real-data-experiment-result-appendix}
\end{table}


\end{document}